\documentclass[10pt]{article}
\usepackage{amsmath, amsthm, amssymb, }
\usepackage{float,epsfig}
\usepackage{color}
\usepackage{subfigure}
\usepackage[utf8]{inputenc}
\usepackage{multirow}
\usepackage{cleveref}

\usepackage{graphicx}
\usepackage{varwidth}
\begin{document}

\newtheorem{theorem}{\bf Theorem}[section]
\newtheorem{proposition}[theorem]{\bf Proposition}
\newtheorem{definition}[theorem]{\bf Definition}
\newtheorem{corollary}[theorem]{\bf Corollary}
\newtheorem{example}[theorem]{\bf Example}
\newtheorem{exam}[theorem]{\bf Example}
\newtheorem{remark}[theorem]{\bf Remark}
\newtheorem{lemma}[theorem]{\bf Lemma}
\newcommand{\nrm}[1]{|\!|\!| {#1} |\!|\!|}

\newcommand{\calL}{{\mathcal L}}
\newcommand{\calX}{{\mathcal X}}
\newcommand{\calA}{{\mathcal A}}
\newcommand{\calB}{{\mathcal B}}
\newcommand{\calC}{{\mathcal C}}
\newcommand{\calK}{{\mathcal K}}
\newcommand{\C}{{\mathbb C}}
\newcommand{\R}{{\mathbb R}}
\renewcommand{\SS}{{\mathbb S}}
\newcommand{\LL}{{\mathbb L}}
\newcommand{\st}{{\star}}
\def\kernel{\mathop{\rm kernel}\nolimits}
\def\sigan{\mathop{\rm span}\nolimits}

\newcommand{\klasse}{{\boldsymbol \Delta}}

\newcommand{\ba}{\begin{array}}
\newcommand{\ea}{\end{array}}
\newcommand{\von}{\vskip 1ex}
\newcommand{\vone}{\vskip 2ex}
\newcommand{\vtwo}{\vskip 4ex}
\newcommand{\dm}[1]{ {\displaystyle{#1} } }

\newcommand{\be}{\begin{equation}}
\newcommand{\ee}{\end{equation}}
\newcommand{\beano}{\begin{eqnarray*}}
\newcommand{\eeano}{\end{eqnarray*}}
\newcommand{\inp}[2]{\langle {#1} ,\,{#2} \rangle}
\def\bmatrix#1{\left[ \begin{matrix} #1 \end{matrix} \right]}
\def \noin{\noindent}
\newcommand{\evenindex}{\Pi_e}



\def \R{{\mathbb R}}
\def \C{{\mathbb C}}
\def \K{{\mathbb K}}
\def \Q{{\mathbb Q}}
\def \H{{\mathbb H}}
\def \L{\mathcal{L}}
\def \pf{{\bf Proof: }}
\def \eproof{$\blacksquare$}
\def \proof{\noin\pf}

\def \T{{\mathbb T}}
\def \Pb{\mathrm{P}}
\def \N{{\mathbb N}}
\def \Ib{\mathrm{I}}
\def \Ls{{\Lambda}_{m-1}}
\def \Gb{\mathrm{G}}
\def \Hb{\mathrm{H}}
\def \Lam{{\Lambda}}

\def \Qb{\mathrm{Q}}
\def \Rb{\mathrm{R}}
\def \Mb{\mathrm{M}}
\def \norm{\nrm{\cdot}\equiv \nrm{\cdot}}

\def \P{{\mathbb{P}}_m(\C^{n\times n})}
\def \A{{{\mathbb P}_1(\C^{n\times n})}}
\def \H{{\mathbb H}}
\def \L{{\mathbb L}}
\def \G{{\F_{\tt{H}}}}
\def \S{\mathbb{S}}
\def \s{\mathbb{s}}
\def \sigmin{\sigma_{\min}}
\def \elam{\Lambda_{\epsilon}}
\def \slam{\Lambda^{\S}_{\epsilon}}
\def \Ib{\mathrm{I}}
\def \Tb{\mathrm{T}}
\def \d{{\delta}}

\def \Lb{\mathrm{L}}
\def \N{{\mathbb N}}
\def \Ls{{\Lambda}_{m-1}}
\def \Gb{\mathrm{G}}
\def \Hb{\mathrm{H}}
\def \Delta{\triangle}
\def \Rar{\Rightarrow}
\def \p{{\mathsf{p}(\lam; v)}}

\def \D{{\mathbb D}}

\def \tr{\mathrm{Tr}}
\def \cond{\mathrm{cond}}
\def \lam{\lambda}
\def \sig{\sigma}
\def \sign{\mathrm{sign}}

\def \ep{\epsilon}
\def \diag{\mathrm{diag}}
\def \rev{\mathrm{rev}}
\def \vec{\mathrm{vec}}

\def \herm{\mathsf{Herm}}
\def \sym{\mathsf{sym}}
\def \odd{\mathsf{sym}}
\def \en{\mathrm{even}}
\def \rank{\mathrm{rank}}
\def \pf{{\bf Proof: }}
\def \dist{\mathrm{dist}}
\def \rar{\rightarrow}

\def \rank{\mathrm{rank}}
\def \pf{{\bf Proof: }}
\def \dist{\mathrm{dist}}
\def \Re{\mathsf{Re}}
\def \Im{\mathsf{Im}}
\def \re{\mathsf{re}}
\def \im{\mathsf{im}}

\def \sym{\mathsf{sym}}
\def \sksym{\mathsf{skew\mbox{-}sym}}
\def \odd{\mathrm{odd}}
\def \even{\mathrm{even}}
\def \herm{\mathsf{Herm}}
\def \skherm{\mathsf{skew\mbox{-}Herm}}
\def \str{\mathrm{ Struct}}
\def \eproof{$\blacksquare$}

\def \bS{{\bf S}}
\def \cA{{\cal A}}
\def \E{{\mathcal E}}
\def \X{{\mathcal X}}
\def \F{{\mathcal F}}
\def \cH{\mathcal{H}}
\def \cJ{\mathcal{J}}
\def \tr{\mathrm{Tr}}
\def \range{\mathrm{Range}}
\def \adj{\star}

\def \pal{\mathrm{palindromic}}
\def \palpen{\mathrm{palindromic~~ pencil}}
\def \palpoly{\mathrm{palindromic~~ polynomial}}
\def \odd{\mathrm{odd}}
\def \even{\mathrm{even}}

\newcommand{\tb}[1]{\textcolor{blue}{ #1}}
\newcommand{\tm}[1]{\textcolor{magenta}{ #1}}
\newcommand{\tre}[1]{\textcolor{red}{ #1}}

\title{Eigenvalue embedding problem for  quadratic regular matrix polynomials with symmetry structures}
 \author{Tinku Ganai\thanks{Department of Mathematics,
IIT Kharagpur, India, E-mail:
tinkuganaimath@gmail.com } \,\, and \,\, Bibhas Adhikari\thanks{Corresponding author, Department of Mathematics, IIT Kharagpur, India, E-mail:
bibhas@maths.iitkgp.ac.in } }

\date{}

\maketitle
\thispagestyle{empty}

{\small \noin{\bf Abstract.} In this paper, we propose a unified approach for solving structure-preserving eigenvalue embedding problem (SEEP) for quadratic regular matrix polynomials  with symmetry structures. First, we determine perturbations of a quadratic matrix polynomial, unstructured or structured , such that the perturbed polynomials reproduce a desired invariant pair while maintaining
the invariance of another invariant pair of the unperturbed polynomial. If the latter is unknown, it is referred to as no
spillover perturbation. Then we use these results for solving the SEEP for structured quadratic matrix polynomials that include: symmetric, Hermitian, $\star$-even and $\star$-odd quadratic matrix polynomials. Finally, we show that the obtained analytical expressions of perturbations can realize existing results for structured polynomials that arise in real-world applications, as special cases. The obtained results are  supported with numerical examples. 


\vone \noin{\bf Keywords.} Eigenvalue embedding problem, model updating, invariant pair, inverse eigenvalue problem, no spillover

\vone\noin{\bf AMS subject classifications.} 15A22, 65F18, 93B55, 47A75

\section{Introduction}

In this paper, we investigate preservation of invariant pairs of unstructured and structured quadratic matrix polynomials under perturbations of the coefficient matrices. This leads to determine perturbations of coefficients of a quadratic matrix polynomial such that a desired set of eigenvalues can be reproduced by the perturbed polynomials that replace a given set of eigenvalues of the unperturbed polynomial. This problem arises in real-world vibrating structural models, and hence we revisit and propose solutions for the well-known quadratic finite element \textit{model updating problem} (MUP) \cite{FrisMotter}. In addition, if the structure-preserving perturbed polynomials preserve the rest of the eigenvalues (need not be known) of the unperturbed structured polynomial then the problem is referred to as \textit{structure-preserving eigenvalue embedding problem} (SEEP)  or MUP with \textit{no spillover} in the literature \cite{Kuo_Datta,ChuAIAA}. We formulate the SEEP in terms of reproducing a desired invariant pair for the perturbed polynomials while preserving another invariant pair of the unperturbed polynomial that need not be known, and consequently we obtain analytical solution for SEEP for a variety of structured quadratic matrix polynomials. Here we mention that SEEP or MUP with no spillover can also be defined as partial inverse eigenvalue problem for structured matrix polynomials \cite{mao2016quadratic}.


Let $\K\in\{\R, \C\},$ where $\R$ and $\C$ denote the field of real and complex numbers respectively. A pair $(X,\Lam)\in \K^{n \times p} \times \K^{p \times p}$ is called an invariant pair \cite{Betcke_Daniel,beyn2010continuation} of  a quadratic matrix polynomial $Q(\lam)=\lam^2 M+\lam D + K\in\K^{n\times n}[\lam]$ if 
\begin{equation} \label{invariant_defn} 
Q(X,\Lam):=MX\Lam^2+DX\Lam+KX=0. 
\end{equation} If $\Lambda = \lam_0 I- J$ where $I$ is the identity matrix of order $p,$ $J$ is the Jordan block of size $p$ and $\lam_0$ is an eigenvalue of $Q(\lam)$ then the columns of $X$ form a Jordan chain of $Q(\lam),$ and the pair $(X, \Lam)$ is called a Jordan pair of $Q(\lam)$ \cite{gohberg2005matrix}. If $p=n$ and $X$ is invertible then $S=X\Lam X^{-1}$ is a solvent, that is, $S$ is a solution of the quadratic matrix equation associated with $Q(\lam)$ \cite{higham2000numerical}. Further, if $p=1$ then $(\Lam,X)$ is an eigenpair of $Q(\lam)$. If $\zeta$ is an eigenvector of $\Lam$ corresponding to the eigenvalue $\lam_0$ then it follows that $(\lam_0,X \zeta)$ is an eigenpair of $Q(\lam)$. Consequently, the eigenvalues of $\Lam$ are eigenvalues of $Q(\lam)$.  Thus invariant pair provides a unified perspective on the problem of computing several eigenpairs for a given matrix polynomial \cite{barkatou2017contour,effenberger2013robust,effenberger2012chebyshev,kressner2009block,ugalde2015computation}. Obviously, invariant pair extends the concepts of standard pair \cite{gohberg2005matrix} and null pair \cite{ball2013interpolation}.  Besides, invariant pairs play a key role in developing several algorithms related to nonlinear eigenvalue problems, for example, see \cite{van2015computing,kressner2014memory}. In \cite{Betcke_Daniel}, the authors analyze the behavior of invariant pairs of a matrix polynomial under the perturbations of the coefficients of the polynomial, and they propose a first-order perturbation expansion. See also 
\cite{szyld2014several}.



 On the other hand, as invariant subspaces of matrices  can be interpreted as a generalization of eigenvectors, deflating pairs and invariant pairs  can be seen as generalizations of eigenpairs for matrix pencils and matrix polynomials respectively. Recently, structure-preserving perturbations for structured matrices and matrix pencils are determined to preserve invariant subspaces and deflating pairs respectively in \cite{ganai2020preserving,Ganai_linearpencil}. Then these results are used to reproduce a desired scalars as eigenvalues of the perturbed structured matrices and matrix pencils without effecting a set of desired eigenpairs of the unperturbed matrix and matrix pencils, respectively. In this paper, we consider the analogous approach for invariant pairs of regular structured matrix polynomials that can be employed to preserve desired spectral properties in the perturbed polynomials. We mention that  Mackey et al. have investigated preserving invariant pairs under Möbius transformation of a matrix polynomial \cite{mackey2015mobius}. Now we describe the structured matrix polynomials which we consider in this paper.

Let $M^*$ and $M^T$ denote the conjugate transpose and transpose of a matrix $M\in \K^{q \times r}$ respectively. Then a polynomial $Q(\lam)=\lam^2 M+\lam D+K \in \K^{n \times n}[\lam]$ is said to have $(\star,\epsilon_1,\epsilon_2)$-structure if 
\begin{equation} \label{Q_nDefn}
M^\star=\epsilon_1 M, \qquad D^\star=\epsilon_2 D, \qquad K^\star=\epsilon_1 K,
\end{equation} 
where $\ep_1, \ep_2\in\{1, -1\},$ $\st\in\{*, T\}.$ These structured quadratic matrix polynomials are known under the following names in the literature, and the symmetry of coefficient matrices induces eigenvalue pairing of the polynomials \cite{mackey2006vector}.
\begin{center}
\begin{tabular}{l|l|l}
name & $(\st,\epsilon_1,\epsilon_2)$ & eigenvalue pairing
\\
\hline
symmetric & $(T,1,1)$ & $\lam$ if $\K=\C$ and $(\lam,\overline{\lam})$ if $\K=\R$ \\
Hermitian & $(*,1,1)$ & $(\lam, \overline{\lam})$ \\
$T$-odd & $(T,-1,1)$ & $(\lam,-\lam)$ if $\K=\C$ and $(\lam,\overline{\lam},-\lam,-\overline{\lam})$ if $\K=\R$\\
$*$-odd & $(*,-1,1)$ & $(\lam,-\overline{\lam})$ \\
$T$-even & $(T,1,-1)$ & $(\lam,-\lam)$ if $\K=\C$ and $(\lam,\overline{\lam},-\lam,-\overline{\lam})$ if $\K=\R$ \\
$*$-even & $(*,1,-1)$ & $(\lam,-\overline{\lam})$ \\
\end{tabular}
\end{center}

Thus we consider the following problem in this paper. \\

\textbf{(P) (Change of invariant pairs with no spillover)} Let $(X_c,\Lam_c)\in \K^{n \times p_1}\times \K^{p_1 \times p_1}$ and $(X_f,\Lam_f)\in \K^{n \times p_2}\times \K^{p_2 \times p_2}$ be invariant pairs of a quadratic matrix polynomial $Q(\lam)=\lam^2 M+\lam D+K \in \K^{n \times n}[\lam]$. Let $(X_a,\Lam_a)$ be a matrix pair of the same dimension as $(X_c,\Lam_c)$. Then find perturbations $\triangle M,\, \triangle D,\, \triangle K \in \K^{n\times n}$ such that $(X_a,\Lam_a)$ and $(X_f,\Lam_f)$ are invariant pairs of a perturbed quadratic matrix polynomial $Q_{\triangle}(\lam)=\lam^2(M+\triangle M)+\lam (D+\triangle D)+(K+\triangle K).$  Moreover, determine perturbations $\triangle M,\,\triangle D,\,\triangle K\in \K^{n\times n}$ such that $Q_{\triangle }(\lam)\in \Q_n(\star,\epsilon_1,\epsilon_2)$ whenever $Q(\lam)\in \Q_n(\star,\epsilon_1,\epsilon_2)$ and $(X_f,\Lam_f)$ need not be known, the phenomena which is known as \textit{no spillover} effect in the literature of structural models. (Note that $\Lam_c,\,\Lam_f,\,\Lam_a$ are not necessarily diagonal. The notations $^c,\,^f,\,^a$ stand for $change,\,fixed$ and $aimed$ respectively.) \\

We call the invariant pairs $(X_c,\Lam_c)$ and $(X_f,\Lam_f)$ of the quadratic matrix polynomial $Q(\lam)$ in \textbf{(P)} as \textit{change} and \textit{fixed} invariant pairs respectively. Note that when $\Lam_c$ and $\Lam_f$ are diagonal matrices, the diagonal entries are eigenvalues of $Q(\lam)$ and the columns of $X_c$ and $X_f$ are eigenvectors associated with those eigenvalues, and hence the Problem \textbf{(P)} boils down to SEEP or MUP with no spillover for $(\st,\ep_1,\ep_2)$-structured quadratic matrix polynomials. The SEEP or MUP arises in quadratic finite element models for real systems in vibration industries  \cite{CarvalhoMSSP, Dattabook, LancasterBook, Mao_Dai, FrisMotter, TisseurQEP}. In particular, MUP has applications in damage detection and health monitoring of vibrating structures such as bridges, highways  \cite{Jaishi, Kim}, and to control resonance vibration of such structures \cite{DattaSarkissian, Elhaybook}. A detailed study on MUP and its applications can be found in the seminal book by Friswell and Mottershead \cite{Mottershedbook}. Indeed, the finite element model corresponding to such a real system is described by a system of second-order ordinary differential equations as
\begin{equation} \label{Quad_diffeqn}
M\ddot{x}(t)+D\dot{x}(t)+Kx(t)=0
\end{equation}    
where $M,\,D,\,K$ are real or complex matrices of order $n\times n$ and $x(t)$ is a column vector of order $n$. Solutions of (\ref{Quad_diffeqn}) can be obtained as $x(t) = x_0e^{\lambda_0t}$, where $(\lam_0, x_0)$ turns out to be eigenpairs of the quadratic matrix polynomial $Q(\lam)=\lam^2 M+\lam D+K$. Further, depending on applications, the matrices $M, D$ and $K$ have certain structures. For example, in vibrating structural systems, $M$ is symmetric positive definite, $K$ is positive semi-definite and $D$ is symmetric and they are called mass, stiffness and damping matrices, respectively \cite{CarvalhoMSSP}. Hence the corresponding $Q(\lam)$ is a symmetric matrix polynomial. On the other hand, for gyroscopic systems, such as rotors of the generator, solar panels on the satellite, $M$ is symmetric positive definite, $K$ is positive semi-definite and $D$ is skew-symmetric known as mass,  stiffness and gyroscopic matrices, respectively \cite{Mao_Dai}, which corresponds to a $T$-even polynomial $Q(\lam)$ \cite{mackey2006vector}. 

There is a voluminous literature on solving MUP or SEEP for vibrating structural models using different methods \cite{Bai_optimization, BaiDattaWang, Baruch, BermanNagy, Brahma_Datta, CarvalhoSoundVibration, Chu_postdefn, chu2007updating, Datta, DattaSarkissianExample, Datta_optimization, Kuo_Datta, Mao_Dai, Zimm_Wind}. However, there are only a few articles that provide analytical expressions of the perturbations $\triangle Q(\lam)=\lam^2\triangle M + \lam\triangle D +\triangle K$ which solve SEEP for structured $Q(\lam)$ \cite{CarvalhoSoundVibration, Chu_postdefn, Kuo_Datta, Mao_Dai}. There are several articles that consider the MUP for symmetric or Hermitian quadratic matrix polynomials \cite{CarvalhoSoundVibration, ChuAIAA, KuoChengLin, lancaster06, lancaster08}. There are a few articles on quadratic MUP with no spillover condition that provide explicit expression of the updating matrices \cite{CarvalhoSoundVibration, Chu_postdefn, Kuo_Datta, Mao_Dai}. However, perturbation matrices presented in these articles can not replace a zero eigenvalue by a nonzero eigenvalue in the perturbed polynomials. For undamped structural models, that is, setting $D=0$ in equation (\ref{Quad_diffeqn}),  solutions of MUP with no spillover effect is obtained recently in  \cite{Ganai_linearpencil} where the corresponding matrix pencils are symmetric, Hermitian, $\st$-even and $\star$-odd. On the other hand, in \cite{Mao_Dai}, the authors consider $T$-even quadratic matrix polynomials with $M$ and $K$ as positive definite or semi-definite matrices respectively.



In this paper, we consider only regular matrix polynomials with the leading coefficient matrix nonsingular. The contribution of our work are as follows. First, analytical expression of perturbations is obtained that solve Problem \textbf{(P)} for unstructured quadratic matrix polynomials when the fixed invariant pair of the polynomial is known. Then we provide explicit expression of structured perturbations that solve the Problem \textbf{(P)} when $Q(\lam) \in \Q_n(\star,\epsilon_1,\epsilon_2)$ and the fixed invariant pair of $Q(\lam)$ need not be known. We utilize these results to obtain parametric solutions for MUP, and SEEP or MUP with no spillover for $(\st,\epsilon_1,\epsilon_2)$-structured quadratic matrix polynomials. Finally, we show that these solutions can identify the solutions for SEEP that are obtained in \cite{Chu_postdefn}, \cite{Kuo_Datta} (see Remark \ref{Reco_Chu_postdefn}), and  in \cite{Mao_Dai} (see Remark \ref{Reco_Mao_Dai}) for specific structured polynomials.

This paper is organized as follows. Next section presents some basic facts on eigenpairs and invariant pairs of quadratic matrix polynomials, and we analyze the connection between invariant pairs and the coefficient matrices of the associated quadratic matrix polynomials having $(\star, \ep_1,\ep_2)$-structure. In Section 3, we provide solution of Problem \textbf{(P)} for unstructured quadratic matrix polynomials when the fixed invariant pair is completely known. In addition, we provide parametric solutions of the Problem \textbf{(P)} for $(\st,\epsilon_1,\epsilon_2)$-structured quadratic matrix polynomials when $(X_f,\Lam_f)$ need not be known. In the next section, we present solutions of the SEEP for specially structured quadratic matrix polynomials. Finally, in Section 5 we illustrate our results with the help of numerical examples. 

\textbf{Notation.} 
We denote $\K^{n \times n}[\lam]$ as the space of one parameter $(\lam)$ matrix polynomials whose coefficient matrices are of order $n \times n$ with its entries are from the field $\K$. We mention that $\mathrm{i}\R$ denote the set of all imaginary numbers. By $\sigma(A)$ we denote the spectrum (set of all eigenvalues) of $A$. $A \geq 0$ denotes that $A$ is a Hermitian positive semi-definite matrix while $A>0$ denote that $A$ is a Hermitian positive definite matrix. $\|X\|_F$ denotes the Frobenius norm of $X$. $\re(x)$ and $\im(x)$ denote the real and imaginary parts of a vector or scalar $x$. Finally, $I_m$ denotes the identity matrix of order $m \times m$.

\section{Invariant pairs of quadratic matrix polynomials}

As mentioned before, we consider regular polynomials $Q(\lam)=\lam^2 M+\lam D+K \in \K^{n \times n}[\lam]$ in this paper, that is, the characteristic polynomial $\chi(\lam)=\det(\lam^2 M+\lam D+K)$ is not a zero polynomial. We also assume that $M$ is nonsingular and hence the roots of the equation $\chi(\lam_0)=0,$ known as eigenvalues of $Q(\lam)$ are all finite. A nonzero vector $x_0\in \K^{n}$ is said to be an eigenvector corresponding to an eigenvalue $\lam_0\in \C$ if $Q(\lam_0,x_0) := (\lam_0^2 M+\lam_0 D+K)x_0=0$. Consequently, $(\lam_0,x_0)$ is called an eigenpair of $Q(\lam)$.  

 It is evident that invariant pair is independent of the choice of the basis: If $(X, \Lam) \in \K^{n \times p} \times \K^{p \times p}$ is an invariant pair of $Q(\lam)$ and $\widetilde{X}=XZ,$ where $Z\in\K^{p\times p}$ is a nonsingular matrix then $(\widetilde{X},\widetilde{\Lam}),$ $\widetilde{\Lam}=Z^{-1}\Lam Z$ is also an invariant pair of $Q(\lam).$ An application of this fact is as follows. Let $Q(\lam)\in \R^{n\times n}[\lam]$ and $(\lam_0,x_0)$ be an eigenpair of $Q(\lam)$ such that $\lam_0\in\C \smallsetminus \R.$ Then  $(\overline{\lam}_0,\overline{x}_0)$ is also an eigenpair of $Q(\lam)$. Now setting $X=\left[x_0\,\, \overline{x}_0 \right]$ and $\Lam=\diag(\lam_0,\overline{\lam}_0)$ it can be checked that $(X,\Lam)\in \C^{n\times 2} \times \C^{2\times 2}$ is an invariant pair of $Q(\lam)$. However, setting $Z=\dfrac{1}{2}\bmatrix{1& -\mathrm{i}\\1& \mathrm{i}},$ a corresponding real  invariant pair is given by $(X_r,\Lam_r)$ of $Q(\lam),$ where $$X_r:=XZ=\bmatrix{\re(x_0)& \im(x_0)}, \,\, \Lam_r:=Z^{-1}\Lam Z=\bmatrix{\re(\lam_0)& \im(\lam_0)\\- \im(\lam_0)& \re(\lam_0)}.$$     

Now we recall the following definition from \cite{Betcke_Daniel}.
\begin{definition}
A pair $(X, \Lam)\in \K^{n\times k}\times \K^{k\times k}$ is called minimal if there exists a positive integer $m$ such that the matrix $$\bmatrix{X\Lam^{m-1} \\ X\Lam^{m-2} \\ \vdots \\ X\Lam \\ X}$$ has full column rank. The smallest such $m$ is called minimality index of $(X, \Lam)$. 
\end{definition}

If $(X_1,\Lam_1)\in\K^{n\times p_1} \times \K^{p_1 \times p_1}$ and $(X_2,\Lam_2)\in\K^{n\times p_2} \times \K^{p_2 \times p_2}$ are two invariant pairs of a polynomial $Q(\lam)=\lam^2 M+\lam D + K\in\K^{n\times n}[\lam]$ then it is easy to verify that $\left(X=[X_1 \,\, X_2], \Lam=\diag(\Lam_1, \Lam_2)\right)\in\K^{n\times (p_1+p_2)}\times \K^{(p_1+p_2)\times (p_1+p_2)}$ is an invariant pair of $Q(\lam).$ It is also known that minimality index of any minimal invariant pair of a quadratic matrix polynomial is less or equal to $2$ [Lemma 5, \cite{Betcke_Daniel}]. Obviously, a sufficient condition for an invariant pair $(X, \Lam)\in\K^{n\times p}\times \K^{p\times p}$ of $Q(\lam)$ is minimal if $X$ is a full column rank matrix. 




\subsection{Structured quadratic matrix polynomials}
Recall that $Q(\lam):=\lam^2 M+ \lam D+K \in \K^{n \times n}[\lam]$ is said to have $(\star,\epsilon_1,\epsilon_2)$-structure if 
\begin{equation*}
M^\star=\epsilon_1 M,\hspace{1cm} D^\star=\epsilon_2 D,\hspace{1cm} K^\star=\epsilon_1 K,  
\end{equation*} 
where $\ep_1, \ep_2\in\{1, -1\}$ and $\st\in\{*, T\}.$ We denote these structured quadratic matrix polynomials as $\Q_n(\star,\epsilon_1,\epsilon_2)$.  Suppose $\lam \in \C$ then we define $\lam^\st=\overline{\lam}$ (the conjugate of $\lam$) if
$\st=*$ and $\lam^\st = \lam$ if $\st = T$. Now we briefly discuss some properties of invariant pairs of structured quadratic matrix polynomials that will be used in sequel.

\begin{proposition} \label{Property}
Let $Q(\lam):=\lam^2 M+\lam D+K\in \Q_n(\star,\epsilon_1,\epsilon_2)$ and $(X_j,\Lam_j)\in \K^{n \times p_j} \times \K^{p_j \times p_j},\,j\in \{1,2\}$ be invariant pairs of $Q(\lam).$ Suppose $S_{jk}:=X_j^\star MX_k \Lam_k+\epsilon_1 \epsilon_2 \Lam_j^\star X_j^\star MX_k+X_j^\star DX_k,\,j,k\in \{1,2\}.$  Then:
\begin{itemize}
\item[$(a)$] $\lam_0$ is an eigenvalue of $Q(\lam)$ if and only if $\ep_1\ep_2 \lam_0^\star$ is an eigenvalue of $Q(\lam)$, 
\item[$(b)$] $S_{jk} \Lam_k=\epsilon_1 \epsilon_2 \Lam_j^\star S_{jk},$
\item[$(c)$] $S_{jk}=0$ whenever $\sigma(\epsilon_1 \epsilon_2 \Lam_j^\star) \cap \sigma(\Lam_k)=\emptyset.$
\end{itemize}
\end{proposition} 
\begin{proof}
Note that $\det\left(\lam_0^2 M+\lam_0 D+K \right)=0 \Leftrightarrow \det \left((\lam_0^\star)^2 M^\star +\lam_0^\star D^\star+K^\star \right)=0 \linebreak \Leftrightarrow \det\left(\epsilon_1(\lam_0^\star)^2 M +\epsilon_2 \lam_0^\star D+\epsilon_1 K \right)=0 \Leftrightarrow \det\left((\epsilon_1 \epsilon_2\lam_0^\star)^2 M +(\epsilon_1 \epsilon_2 \lam_0^\star) D+K \right)=0$. Thus $(a)$ follows. As $(X_j,\Lam_j)$ is an invariant pair of $Q(\lam)$ so we have $Q(X_j,\Lam_j)=MX_j\Lam_j^2+DX_j\Lam_j+KX_j=0,$ that is $-KX_j=MX_j\Lam_j^2+DX_j\Lam_j$ then operating $\star$ on it we get $-\epsilon_1X_j^\star K=\epsilon_1 (\Lam_j^\star)^2 X_j^\star M+\epsilon_2 \Lambda_j^\star X_j^\star D$ and postmultiplying it by $X_k$ we obtain 
\begin{equation}\label{ortho1}
-X_j^\star KX_k=(\Lam_j^\star)^2 X_j^\star MX_k+\epsilon_1 \epsilon_2 \Lambda_j^\star X_j^\star DX_k.
\end{equation}
Since $(X_k,\Lam_k)$ is an invariant pair of $Q(\lam)$ so we have $Q(X_k,\Lam_k)=MX_k\Lam_k^2+DX_k\Lam_k+KX_k=0$ then premultiplying it by $X_j^\star$ it gives 
\begin{equation} \label{ortho2}
-X_j^\star KX_k= X_j^\star MX_k \Lam_k^2+X_j^\star DX_k \Lam_k.
\end{equation}
From $(\ref{ortho1})$ and $(\ref{ortho2})$ it follows that 
\begin{equation} \label{ortho3}
(\Lam_j^\star)^2 X_j^\star MX_k+\epsilon_1 \epsilon_2 \Lambda_j^\star X_j^\star DX_k=X_j^\star MX_k \Lam_k^2+X_j^\star DX_k \Lam_k.\end{equation}
Then adding both sides of $(\ref{ortho3})$ by $\epsilon_1 \epsilon_2 \Lam_j^\star X_j^\star MX_k \Lam_k$ we get 
$$\epsilon_1 \epsilon_2 \Lam_j^\star \left(X_j^\star MX_k \Lam_k+\epsilon_1 \epsilon_2 \Lam_j^\star X_j^\star MX_k+X_j^\star DX_k \right)= \left(X_j^\star MX_k \Lam_k+\epsilon_1 \epsilon_2 \Lam_j^\star X_j^\star MX_k+X_j^\star DX_k \right) \Lam_k$$
that is $S_{jk} \Lam_k=\epsilon_1 \epsilon_2 \Lam_j^\star S_{jk}$ hence $(b)$ follows. Now solving the homogeneous Sylvester equation   $\epsilon_1 \epsilon_2 \Lam_j^\star S_{jk}-S_{jk} \Lam_k=0$ we get $S_{jk}=0$ whenever $\sigma(\epsilon_1 \epsilon_2 \Lam_j^\star) \cap \sigma(\Lam_k)= \emptyset.$ Hence $(c)$ follows. $\hfill{\square}$
\end{proof}
Then we have the following corollary.
\begin{corollary}\label{S_jProperty}
From Proposition \ref{Property} it follows that $S_{jj}=X_j^\star MX_j \Lam_j+\epsilon_1 \epsilon_2 \Lam_j^\star X_j^\star MX_j+X_j^\star DX_j$ for $j\in \{1,2\}$. Then $S_{jj}^\star= \epsilon_2 X_j^\star MX_j \Lam_j+\epsilon_1 \Lam_j^\star X_j^\star MX_j+\epsilon_2 X_j^\star DX_j=\epsilon_2 S_{jj}.$ Further by Proposition \ref{Property} $(b)$ it follows that $S_{jj}\Lam_j= \epsilon_1 \epsilon_2 \Lam_j^\star S_{jj}=\epsilon_1 (S_{jj} \Lam_j)^\star$.
\end{corollary} 


Suppose $\lam\in \C$ and $x\in \C^{n}$ then using Corollary \ref{S_jProperty} we note that 
$$(\lam + \epsilon_1 \epsilon_2 \lam^\star) x^\star Mx+x^\star Dx \begin{cases} \in \R &\text{if }(\st, \epsilon_2)=(*,1),\\
 \in \mathrm{i}\R &\text{if }(\st, \epsilon_2)=(*,-1),\\
 =0&\text{if }(\st, \epsilon_2)=(T,-1).
\end{cases}
$$
It should be noted that the matrices $X_1$ and $X_2$ in Proposition \ref{Property} may be identical, then we have the following corollary.
\begin{corollary}
Let $(X,\Lam)$ be an invariant pair of $\lam^2 M+\lam D+K \in \Q_n(\star,\epsilon_1,\epsilon_2)$ with $\sigma(\epsilon_1 \epsilon_2 \Lam^\star) \cap \sigma(\Lam)=\emptyset$ leads to $X^\star MX \Lam+\epsilon_1 \epsilon_2 \Lam^\star X^\star MX+X^\star DX=0.$
\end{corollary}
On choosing $X_1,\,X_2$ in Proposition \ref{Property}
as column vectors we have the following corollary.
\begin{corollary} \label{ortho_eig}
Let $(\lam_1,x_1)$ and $(\lam_2,x_2)$ be eigenpairs of $\lam^2 M+\lam D+K \in \Q_n(\star,\epsilon_1,\epsilon_2).$ Then $\lam_2 x_1^\star Mx_2 +\epsilon_1 \epsilon_2 \lam_1^\star x_1^\star Mx_2+x_1^\star Dx_2=0$ whenever $\lam_2 \neq \epsilon_1 \epsilon_2 \lam^\star_1$.
\end{corollary} 


\begin{corollary} \label{blockform}
Let $(\lam_0,x_0)$ be an eigenpair of $Q(\lam):=\lam^2 M+\lam D+K \in \Q_n(\star,\epsilon_1,\epsilon_2)$ with $\lam_0\neq \epsilon_1 \epsilon_2 \lam_0^\star$. Then Proposition \ref{Property} $(a)$ implies that $\epsilon_1 \epsilon_2 \lam_0^\star$ is also an eigenvalue of $Q(\lam)$ with corresponding eigenvector $\tilde{x}_0$, that is $(\epsilon_1 \epsilon_2 \lam_0^\star,\tilde{x}_0)$ is an eigenpair of $Q(\lam).$ Clearly $(X_0,\Lam_0)$ is an invariant pair of $Q(\lam)$ where $X_0:=\bmatrix{x_0& \tilde{x}_0}$ and $\Lam_0:=\diag(\lam_0,\epsilon_1 \epsilon_2 \lam_0^\star)$. Setting $s_0:=2 \epsilon_1 \epsilon_2 \lam_0^\star x_0^\star M \tilde{x}_0+x_0^\star D \tilde{x}_0$ we obtain,
\begin{center}
$X_0^\star MX_0 \Lam_0+\epsilon_1 \epsilon_2 \Lam_0^\star X_0^\star MX_0+X_0^\star DX_0=\bmatrix{0& s_0\\ \epsilon_2 s_0^\star& 0}$.
\end{center}
\end{corollary}
  
\section{Preserving invariant pairs under perturbations}

In this section we determine perturbations of unstructured and structured quadratic matrix polynomials that reproduce a desired matrix pair as an invariant pair of the perturbed polynomials and preserve a desired invariant pair of the corresponding unperturbed polynomial. 

\subsection{Unstructured perturbations}
First we show that given any minimal pair $(X, \Lam)\in \K^{n\times p}\times \K^{p\times p}$ can be reproduced as an invariant pair of infinitely many matrix polynomials as follows.

\begin{proposition}
Let $(X, \Lam)\in \K^{n\times p}\times \K^{p\times p}$ be a minimal pair. Then any polynomial $Q(\lam)=\lam^2 M + \lam D +K\in \K^{n\times n}[\lam]$ such that $(X, \Lam)$ can be an invariant pair of $Q(\lam)$ is given by 
\beano
M &=& Z_1\left(I_n-(X\Lam^2)Q_X(X\Lam^2)^* \right) - Z_2(X\Lam) Q_X(X\Lam^2)^* - Z_3XQ_X(X\Lam^2)^* ,  \\
D &=& -Z_1(X\Lam^2)Q_X(X\Lam)^* + Z_2\left(I_n-(X\Lam) Q_X(X\Lam)^* \right) - Z_3X Q_X(X\Lam)^* , \\
K &=& -Z_1 (X\Lam^2)Q_XX^* - Z_2(X\Lam)Q_XX^* + Z_3(I_n-XQ_XX^*),
\eeano where $Q_X= [(X\Lam^2)^*(X\Lam^2) +(X\Lam)^*(X\Lam) + X^* X]^{-1}$ and  $Z_j\in\K^{n\times n},$ $j=1,\,2,\,3$.
\end{proposition}

\pf The pair $(X,\Lam)$ is an invariant pair of some $Q(z)=z^2 M + z D +K$ if and only if $$ Q(X,\Lam) =0 \,\, \Leftrightarrow \,\, \bmatrix{M & D & K} \bmatrix{X\Lam^2 \\ X\Lam \\ X}=0,$$ which is a homogeneous linear system $AY=0.$ Any solution of this system is of the form $A=Z(I_{3n}-YY^\dagger),$ where $Y^\dagger=(Y^* Y)^{-1}Y^*$ is the pseudoinverse of $Y,$ and $Z\in\K^{n\times 3n}$ is arbitrary. Thus $$\bmatrix{M & D & K} = Z\left( I_{3n} - \bmatrix{(X\Lam^2) Q_X (X\Lam^2)^* & (X\Lam^2)Q_X (X\Lam)^* & (X\Lam^2)Q_XX^* \\ (X\Lam) Q_X (X\Lam^2)^* & (X\Lam) Q_X (X\Lam)^* & (X\Lam) Q_XX^* \\ X Q_X (X\Lam^2)^* & X Q_X (X\Lam)^* & XQ_XX^*}\right),$$ where $Q_X= [(X\Lam^2)^*(X\Lam^2) +(X\Lam)^*(X\Lam) + X^* X]^{-1}.$ Then the desired result follows by writing $Z=[Z_1 \,\, Z_2 \,\, Z_3],$ $Z_j\in\K^{n\times n}.$ \hfill{$\square$}

Recall that the MUP or SEEP is concerned with finding perturbations of coefficients of a given quadratic matrix polynomial $Q(\lam)$ such that a set of known eigenvalues $\lam_i^c,$ $i=1,\hdots,p$ of $Q(\lam)$ are changed by a desired set of compatible scalars $\lam_i^a, i=1,\hdots,p$ as eigenvalues of perturbed quadratic matrix polynomials $Q_\triangle(\lam)$. If $x_i^c$ is an eigenvector corresponding to the eigenvalue $\lam_i^c$ of $Q(\lam)$ then setting $\Lam_c=\diag(\lam_1^c,\hdots,\lam_p^c)$ and $X_c=\bmatrix{x_1^c & \hdots & x_p^c}$ the MUP can be stated as changing the invariant pair $(X_c,\Lam_c)$ of $Q(\lam)$ by $(X_c,\Lam_a)$ as invariant pair of the perturbed polynomials $Q_\triangle(\lam)$, where $\Lam_a=\diag(\lam_1^a,\hdots,\lam_p^a).$ The following theorem provides explicit analytical expression of perturbations of the coefficients of $Q(\lam)$ that solves this problem when $\Lam_c, \Lam_a$ are not necessarily diagonal matrices.


\begin{theorem}
Let $(X_c,\Lam_c)\in \K^{n \times p} \times \K^{p \times p}$ be an invariant pair of the matrix polynomial $Q(\lam)=\lam^2 M+\lam D+K \in \K^{n \times n}[\lam]$. Let $\Lam_a\in \K^{p \times p}$ be such that $(X_c,\Lam_a)$ is minimal. Then any matrix polynomial $Q_\triangle (\lam)=\lam^2(M+\triangle M)+\lam (D+\triangle D)+(K+\triangle K) \in \K^{n \times n}[\lam]$ for which $(X_c,\Lam_a)$ is an invariant pair is given by
$$\triangle M=Z_1-WR(\Lam_a^2)^*X_c^*,\,\,\,\,\triangle D=Z_2-WR \Lam_a^*X_c^*,\,\,\,\,\triangle K=Z_3-WRX_c^*,$$
where
$W=MX_c(\Lam_a^2-\Lam_c^2)+DX_c(\Lam_a-\Lam_c)+Z_1X_c\Lam_a^2+Z_2X_c \Lam_a+Z_3 X_c,\,R=((X_c \Lam_a^2)^*X_c \Lam_a^2 + (X_c \Lam_a)^*X_c \Lam_a + X_c^*X_c)^{-1}$ and $Z_j \in \K^{n \times n},\,j=1,\,2,\,3$ are arbitrary.
\end{theorem}

\pf Since $(X_c,\Lam_c)$ is an invariant pair of $Q(\lam)$, that is, $Q(X_c,\Lam_c)=0,$ hence $KX_c=-MX_c\Lam_c^2-DX_c \Lam_c$. Then $(X_c,\Lam_a)$ is an invariant pair of $Q_\triangle (\lam)$ if and only if the matrices $\triangle M,\,\triangle D,\,\triangle K$ satisfy
\begin{equation} \label{MUP_unstructured_eqn}
\underbrace{\bmatrix{\triangle M & \triangle D & \triangle K}}_{A} \underbrace{ \bmatrix{X_c \Lam_a^2 \\ X_c \Lam_a \\X_c}}_{X}=\underbrace{ MX_c(\Lam_c^2-\Lam_a^2)+DX_c(\Lam_c-\Lam_a)}_{B}.    
\end{equation}

Then the equation $(\ref{MUP_unstructured_eqn})$ has a solution if and only if it satisfies $BX^{\dagger} X=B$ and any solution can be written as  \begin{equation}\label{eqn:map} A=BX^{\dagger}+Z(I_{3n}-XX^{\dagger}),\end{equation} where $Z\in \K^{n \times 3n}$ is arbitrary and
$X^{\dagger}$ denotes the Moore-Penrose pseudoinverse of $X$  \cite{Sun}. As the matrix pair $(X_c,\Lam_a)$ is minimal so its minimality index must be less than or equal to $2,$ so it implies that $X$ is a full column rank matrix, and hence $X^{\dagger}=(X^* X)^{-1} X^*=RX^*$, thus $BX^{\dagger} X=B$ holds. Therefore the desired result follows by setting $Z=[Z_1 \,\,Z_2\,\,Z_3],\,Z_j \in \K^{n \times n}$. $\hfill{\square}$


In the following theorem we determine perturbations of a quadratic matrix polynomial such that the perturbed polynomials change an invariant pair of the corresponding unperturbed polynomial by a desired invariant pair while preserving another invariant pair of the unperturbed polynomial. Thus the following theorem provides solution for the Problem \textbf{(P)} for unstructured quadratic matrix polynomials.

\begin{theorem}
Let $(X_c,\Lam_c) \in \K^{n \times p_1} \times \K^{p_1 \times p_1}$ and $(X_f,\Lam_f) \in \K^{n \times p_2} \times \K^{p_2 \times p_2}$ be two invariant pairs of $Q(\lam)=\lam^2 M+\lam D+K \in \K^{n \times n}[\lam]$. Let $(X_a,\Lam_a) \in \K^{n \times p_1} \times \K^{p_1 \times p_1}$ be a matrix pair such that $(X=[X_a \,\,X_f], \Lam = \mathrm{diag}(\Lam_a,\Lam_f))$ is minimal. Then the perturbed polynomials $Q_{\triangle}(\lam)=\lam^2(M+ \triangle M)+\lam(D+\triangle D)+(K+\triangle K)\in \K^{n \times n}[\lam]$ that reproduce the pairs  $(X_a,\Lam_a)$ and $(X_f, \Lam_f)$ as invariant pairs, are given by 
$$ \triangle M = Z_1 - F_M, \,\, \triangle D = Z_2 - F_D, \,\, \triangle K = Z_3 - F_K,$$
where
\beano
F_M &=& (Q(X_a,\Lam_a) + Z\widetilde{X}_a)(U(X_a \Lam_a^2)^*+V(X_f \Lam_f^2)^*) + (Z \widetilde{X}_f ) (V^* (X_a \Lam_a^2)^*+W(X_f \Lam_f^2)^*), \\
F_D &=& (Q(X_a,\Lam_a)+ Z \widetilde{X}_a)(U(X_a \Lam_a)^*+V(X_f \Lam_f)^*) + (Z \widetilde{X}_f) (V^* (X_a \Lam_a)^*+W(X_f \Lam_f)^*),\\
F_K &=& (Q(X_a,\Lam_a) +Z\widetilde{X}_a)(UX_a^*+VX_f^*) + (Z \widetilde{X}_f) (V^* X_a^*+WX_f^*),
\eeano with $Z=[Z_1 \,\, Z_2 \,\, Z_3], Z_j \in \K^{n \times n},\,j=1,\,2,\,3$ are arbitrary, $$\widetilde{X}_a = \bmatrix{X_a \Lam_a^2 \\ X_a \Lam_a \\ X_a}, \widetilde{X}_f = \bmatrix{X_f \Lam_f^2 \\ X_f \Lam_f \\ X_f}, \widetilde{X}=\bmatrix{X \Lam^2 \\ X \Lam \\ X},\,\left[
\begin{array}{c|c}
U_{p_1 \times p_1} & V_{p_1 \times p_2} \\
\hline
V^*_{p_2 \times p_1} & W_{p_2 \times p_2}
\end{array}
\right] =(\widetilde{X}^* \widetilde{X})^{-1}.$$



\end{theorem}

\pf
Since $(X_c, \Lam_c)$ and $(X_f, \Lam_f)$ are invariant pairs of $Q(\lam)$ so we have \begin{equation} \label{uns_given_eqn}
Q(X_c,\Lam_c)=0= Q(X_f,\Lam_f).
\end{equation}
As $(X_a,\Lam_a)$ and $(X_f,\Lam_f)$ need to be invariant pairs of the updated matrix polynomial $Q_{\triangle}(\lam)=\lam^2(M+\triangle M)+\lam(D+\triangle D)+(K+\triangle K)$ so the matrices $\triangle M,\,\triangle D,\,\triangle K$ must satisfy 
\begin{equation}\label{uns_needed_eqn1}
Q_\triangle(X_a,\Lam_a)=0= Q_\triangle(X_f,\Lam_f).
\end{equation}

Then by equation $(\ref{uns_given_eqn})$  
\begin{equation}\label{uns_needed_eqn2}
\triangle Q(X_a,\Lam_a) = -Q(X_a,\Lam_a), \,\,\, \triangle Q(X_f,\Lam_f)=0,
\end{equation}
where $\triangle Q(\lam)=\lam^2\triangle M + \lam \triangle D +\triangle K$. Thus equation $(\ref{uns_needed_eqn2})$ can be written as,
\begin{equation}\label{uns_sol1}
\underbrace{\bmatrix{\triangle M & \triangle D& \triangle K}}_{A}
\widetilde{X} =
\underbrace{\bmatrix{-Q(X_a,\Lam_a)& 0}}_{B}.
\end{equation}

Then the desired result follows by equation (\ref{eqn:map}). \hfill{$\square$}

It can be seen that the unstructured solution matrices $\triangle M,\,\triangle D,\,\triangle K$ require the knowledge of $(X_f,\Lam_f).$ However, this information is not always available in real world applications. In the next subsection on $(\st,\epsilon_1,\epsilon_2)$-structured quadratic matrix polynomials we find structured perturbations whose construction requires only the knowledge of invariant pair $(X_c,\Lam_c)$ and a property of spectrum $\sigma(\Lam_f)$ which is generically satisfied.


\subsection{Structured perturbations}
First we consider the problem of determining $(\st,\epsilon_1,\epsilon_2)$-structured polynomials that preserve a desired a pair $(X,\Lam)\in\K^{n\times p}\times \K^{p\times p}, p\leq n$ as invariant pair. 
 
\begin{theorem} \label{Str_Inv}
Let $(X,\Lam)\in\K^{n\times p}\times \K^{p\times p},\, p\leq n$ and $\rank(X)=p.$ Then $X$ can be factorized as $X=\bmatrix{Q_1& Q_2} \bmatrix{R\\0}$ where $Q=\bmatrix{Q_1& Q_2} \in \K^{n \times n}$ satisfies $QQ^\st=Q^\st Q=I_n$ and $R\in \K^{p \times p}$ is nonsingular. Then a set of polynomials $Q(\lam)=\lam^2 M+\lam D + K\in\Q_n(\st,\epsilon_1,\epsilon_2)$ such that $(X,\Lam)$ is an invariant pair of $Q(\lam)$ is given by $$M=Q \bmatrix{0 & \epsilon_1M_{12}^\st \\ M_{12} & M_{22}} Q^\st,\,\,\,\,\,\, D=Q \bmatrix{0 & \epsilon_2 D_{12}^\st \\ D_{12} & D_{22}} Q^\st,\,\,\,\,\,\, K=Q \bmatrix{0 & \epsilon_1K_{12}^\st \\ K_{12} & K_{22}} Q^\st$$ where \beano M_{12} &=& [Z_1(I_p-\Lam^2S(\Lam^2)^*) - Z_2 \Lam S(\Lam^2)^* - Z_3 S(\Lam^2)^*]R^{-1}, \\ D_{12} &=& [-Z_1\Lam^2S \Lam^* + Z_2 (I_p-\Lam S \Lam^*) - Z_3 S \Lam^*]R^{-1}, \\ 
K_{12} &=& [-Z_1\Lam^2 S - Z_2 \Lam S + Z_3(I_p-S)]R^{-1},\eeano  $M_{22}^\st=\epsilon_1 M_{22},$ $D_{22}^\st=\epsilon_2 D_{22},$ $K_{22}^\st=\epsilon_1 K_{22}$ are arbitrary matrices of order $(n-p)\times (n-p),$ and $Z_1, Z_2, Z_3\in\K^{(n-p) \times p},\,S=\left(( \Lam^2)^* \Lam^2+\Lam^* \Lam+I_p \right)^{-1}$.
\end{theorem}

\pf As $\rank(X)=p$ so there exists a nonsingular matrix $R\in \K^{p \times p}$ such that $X=\bmatrix{Q_1& Q_2} \bmatrix{R\\ 0}$ holds where $Q=\bmatrix{Q_1& Q_2} \in \K^{n \times n}$ satisfies $Q^\st Q=I_n=QQ^\st$.  
Then the pair $(X,\Lam)$ is an invariant pair of $Q(\lam)$ if and only if \begin{eqnarray} && Q^{\st} MQQ^\st X\Lam^2 + Q^{\st} DQ Q^\st X\Lam + Q^{\st} KQ Q^\st X =0 \nonumber \\ &\Rightarrow& \bmatrix{M_{11} R\Lam^2 \\ M_{12} R\Lam^2} + \bmatrix{D_{11} R\Lam \\ D_{12} R \Lam} + \bmatrix{K_{11} R\\ K_{12} R} =\bmatrix{0 \\ 0}\label{eqn:strm}, 
\end{eqnarray} where $$M=Q \bmatrix{M_{11} & \epsilon_1M_{12}^\st \\ M_{12} & M_{22}} Q^{\star},\, D= Q \bmatrix{D_{11} & \epsilon_2D_{12}^\st \\ D_{12} & D_{22}}Q^{\star},\, K=Q \bmatrix{K_{11} & \epsilon_1K_{12}^\st \\ K_{12} & K_{22}}Q^{\star}.$$ Then setting $M_{11}=D_{11}=K_{11}=0,$ the equation (\ref{eqn:strm}) reduces to solving $$\bmatrix{M_{12} & D_{12} & K_{12}} R \bmatrix{\Lam^2 \\ \Lam \\ I_p}=0.$$ Then by equation (\ref{eqn:map}) the desired result follows. \hfill{$\square$}

Now we determine structure-preserving perturbations of a polynomial $Q(\lam)\in\Q_n(\st,\epsilon_1,\epsilon_2)$ such that the perturbed polynomials $Q_\triangle(\lam)\in\Q_n(\st,\epsilon_1,\epsilon_2)$ replace a known invariant pair $(X_c,\Lam_c)$ of $Q(\lam)$ and preserve a desired invariant pair $(X_c,\Lam_a)$ when $\Lam_a$ has the same dimension as $\Lam_c$ and $X_c$ is a full column rank matrix. 


\begin{theorem}\label{thm:mupip}
Let $(X_c,\Lam_c)\in \K^{n \times p} \times \K^{p \times p}$ be an invariant pair of the polynomial $Q(\lam)=\lam^2 M+\lam D+K \in \Q_n(\st,\ep_1,\ep_2)$ with $\rank(X_c)=p \leq n$. Suppose $X_c=\bmatrix{Q_1& Q_2} \bmatrix{R\\ 0}$ where $R\in \K^{p\times p}$ is nonsingular and $Q=\bmatrix{Q_1& Q_2} \in\K^{n\times n}$ satisfies $Q^\st Q=QQ^\st=I_n$. Let $\Lam_a \in \K^{p \times p}$. Set $$
\triangle M=Q \bmatrix{-Q_1^\st M Q_1& \ep_1 M_{12}^\star \\ M_{12}& M_{22}}Q^\st,\triangle D=Q \bmatrix{-Q_1^\st D Q_1& \ep_2 D_{12}^\star \\ D_{12}& D_{22}}Q^\st,\triangle K=Q \bmatrix{-Q_1^\st K Q_1& \ep_1 K_{12}^\star \\ K_{12}& K_{22}}Q^\st$$
with $$M_{12}=[Z_1-WS(\Lam_a^2)^*]R^{-1},\,\,\,\,D_{12}=[Z_2-WS \Lam_a^*]R^{-1},\,\,\,\,K_{12}=[Z_3-WS]R^{-1},$$
$M_{22}=\ep_1 M_{22}^\star,\,D_{22}=\ep_2 D_{22}^\star,\,K_{22}=\ep_1 K_{22}^\star$ are arbitrary matrices of order $(n-p) \times (n-p)$ and $W=Q_2^\st MX_c(\Lam_a^2-\Lam_c^2)+Q_2^\st DX_c(\Lam_a-\Lam_c)+Z_1\Lam_a^2+Z_2\Lam_a+Z_3,\,S=((\Lam_a^2)^*\Lam_a^2+\Lam_a^* \Lam_a+I_p)^{-1},\,Z_1,\,Z_2,\,Z_3 \in \K^{(n-p) \times p}$. Then $(X_c,\Lam_a)$ is an invariant pair of $Q_\triangle (\lam)=\lam^2 (M+\triangle M)+\lam (D+\triangle D)+ (K+\triangle K) \in \Q_n(\st,\ep_1,\ep_2)$.  
\end{theorem}

\pf Note that $KX_c=-MX_c \Lam_c^2-DX_c\Lam_c$ since $(X_c,\Lam_c)$ is an invariant pair of $Q(\lam),$ that is, $Q(X_c,\Lam_c)=0.$ Now $(X_c,\Lam_a)$ is an invariant pair of $Q_\triangle (\lam)=\lam^2 (M+\triangle M)+\lam (D+\triangle D)+ (K+\triangle K)$ if and only if the matrices $\triangle M,\,\triangle D,\,\triangle K$ satisfy 
\begin{eqnarray} 
    Q^\st \triangle M QQ^\st X_c \Lam_a^2+Q^\st \triangle D QQ^\star X_c \Lam_a+Q^\st \triangle K QQ^\star X_c=Q^\st[MX_c(\Lam_c^2-\Lam_a^2)+DX_c (\Lam_c-\Lam_a)] \nonumber \\ \Rightarrow \bmatrix{M_{11}R \Lam_a^2\\M_{12}R \Lam_a^2}+\bmatrix{D_{11}R \Lam_a \\D_{12}R \Lam_a}+\bmatrix{K_{11}R\\K_{12}R}=\bmatrix{Q_1^\st MQ_1R(\Lam_c^2-\Lam_a^2)+Q_1^\st DQ_1R(\Lam_c-\Lam_a)\\ Q_2^\st \left(MX_c(\Lam_c^2-\Lam_a^2)+DX_c(\Lam_c-\Lam_a)\right)} \label{eqn:mupip}
\end{eqnarray}
where $\triangle M=Q\bmatrix{M_{11}& \ep_1 M_{12}^\st \\ M_{12}& M_{22}}Q^\st,\,\triangle D=Q\bmatrix{D_{11}& \ep_2 D_{12}^\st \\ D_{12}& D_{22}}Q^\st,\,\triangle K=Q\bmatrix{K_{11}& \ep_1 K_{12}^\st \\ K_{12}& K_{22}}Q^\st$. Since $X_c=\bmatrix{Q_1& Q_2}\bmatrix{R\\0}=Q_1R$ and $KX_c=-MX_c \Lam_c^2-DX_c\Lam_c$ so setting $M_{11}=-Q_1^\st MQ_1,\,D_{11}=-Q_1^\st DQ_1,\,K_{11}=-Q_1^\st KQ_1$ the equation (\ref{eqn:mupip}) boils down to solving
$$\bmatrix{M_{12}& D_{12}& K_{12}} R \bmatrix{\Lam_a^2 \\ \Lam_a \\ I_p}=Q_2^\star\left(MX_c(\Lam_c^2-\Lam_a^2)+DX_c(\Lam_c-\Lam_a) \right).$$ Hence the desired result follows by using equation $(\ref{eqn:map})$. $\hfill{\square}$ 

Now we present the solution of Problem \textbf{(P)} for $(\star,\epsilon_1,\epsilon_2)$-structured quadratic matrix polynomials when $(X_f,\Lam_f)$ need not be known. When $X_f$ is completely unknown, we assume that $X_a=X_cP$ that is, $\mathfrak{R}(X_a)=\mathfrak{R}(X_c)$ where $\mathfrak{R}(X)$ denotes the range space of $X$ and $P\in \K^{p_1 \times p_1}$ is nonsingular. First, we consider $P=I_{p_1},$ that is, in the perturbed polynomials an invariant pair $(X_c, \Lam_c)$ is replaced by $(X_c,\Lam_a),$ $\Lam_a\neq \Lam_c$ while keeping another invariant pair $(X_f,\Lam_f)$ of the unperturbed polynomial unchanged. 

\begin{theorem} \label{P2_soln}
Suppose $(X_c,\Lam_c)\in\K^{n\times p_1}\times \K^{p_1\times p_1}$ and $(X_f,\Lam_f)\in\K^{n\times p_2} \times \K^{p_2\times p_2}$ are invariant pairs of the quadratic matrix polynomial $Q(\lam)=\lam^2 M+\lam D+K \in \Q_n(\star,\epsilon_1,\epsilon_2).$ Let  $\Lam_a \in \K^{p_1 \times p_1},$ and denote $S =X_c^\star MX_c \Lam_c+\epsilon_1 \epsilon_2 \Lam_c^\star X_c^\star MX_c+X_c^\star DX_c$. Assume that \begin{enumerate}
    \item $\sig(\Lam_c)\cap \sigma(\epsilon_1 \epsilon_2 \Lam_f^\star) =\emptyset,$
    \item $R:=X_c^\star MX_c \Lam_a+\epsilon_1 \epsilon_2 \Lam_c^\star X_c^\star MX_c+X_c^\star DX_c$ is nonsingular.
\end{enumerate}


Then $(X_c,\Lam_a)$ and $(X_f,\Lam_f)$ are invariant pairs of the quadratic matrix polynomial $Q_{\triangle}(\lam)=\lam^2 (M+\triangle M)+\lam (D+\triangle D)+(K+\triangle K)\in \K^{n\times n}[\lam]$ where
\beano
\triangle M &=& MX_c Z X_c^\star M, \\
\triangle D &=& \epsilon_1 \epsilon_2 MX_c Z \Lam_c^\star X_c^\star M+MX_c Z X_c^\star D+MX_c \Lam_c Z X_c^\star M+DX_c Z X_c^\star M, \\
\triangle K &=& \epsilon_1 \epsilon_2 MX_c \Lam_c Z \Lam_c^\star X_c^\star M+MX_c \Lam_c Z X_c^\star D+\epsilon_1 \epsilon_2 DX_c Z \Lam_c^\star X_c^\star M+DX_c Z X_c^\star D,
\eeano $Z=(\Lam_c-\Lam_a)R^{-1}.$  Moreover, if $S\Lam_a=\epsilon_1 (S\Lam_a)^\star$ then $Q_{\triangle}(\lam) \in \Q_n(\star,\epsilon_1,\epsilon_2)$.


\end{theorem}
\begin{proof}
If $(X_c,\Lam_c)$ and $(X_f,\Lam_f)$ are invariant pairs of $Q(\lam)$ then 
\begin{equation} \label{given_eqn}
Q(X_c,\Lam_c) =0= Q(X_f,\Lam_f).
\end{equation}
Consequently, by Proposition \ref{Property} $(c)$ we have $X_c^\star MX_f \Lam_f+\epsilon_1 \epsilon_2 \Lam_c^\star X_c^\star MX_f+X_c^\star DX_f=0$ whenever $\sigma(\Lam_c) \cap \sigma(\epsilon_1 \epsilon_2 \Lam_f^\star)=\emptyset$. 
Further, $(X_f,\Lam_f)$ is an invariant pair of $Q_{\triangle}(\lam)$ if $Q_\triangle(X_f,\Lam_f)=0,$ which means that $\triangle M,\,\triangle D,\,\triangle K$ must  satisfy $\triangle Q(X_f,\Lam_f) =0,$ by equation $(\ref{given_eqn})$ where $\triangle Q(\lam)= \lam^2 \triangle M+\lam \triangle D+\triangle K$. 

Setting the matrices $\triangle M,\,\triangle D,\,\triangle K$ as described in the statement of the theorem, we obtain 
\beano
&& \triangle Q(X_f,\Lam_f) \\
&=& MX_c Z\left(X_c^\star MX_f \Lam_f+\epsilon_1 \epsilon_2 \Lam_c^\star X_c^\star MX_f+X_c^\star DX_f \right) \Lam_f + MX_c\Lam_c Z (X_c^\star MX_f \Lam_f \\ && +\epsilon_1 \epsilon_2 \Lam_c^\star X_c^\star MX_f+X_c^\star DX_f ) +DX_c Z \left(X_c^\star MX_f \Lam_f+\epsilon_1 \epsilon_2 \Lam_c^\star X_c^\star MX_f+X_c^\star DX_f \right) \\
&=& 0,
\eeano 
where the last equality follows by using $X_c^\star MX_f \Lam_f+\epsilon_1 \epsilon_2 \Lam_c^\star X_c^\star MX_f+X_c^\star DX_f=0$. Thus $(X_f,\Lam_f)$ is an invariant pair of $Q_{\triangle }(\lam)=\lam^2 (M+\triangle M)+\lam (D+\triangle D)+(K+\triangle K)$.


Next, using $KX_c=-MX_c \Lam_c^2-DX_c \Lam_c$ we obtain
\beano
&& (M+\triangle M)X_c\Lam_a^2+(D+\triangle D)X_c \Lam_a+ (K+\triangle K)X_c \\
&=& \triangle MX_c \Lam_a^2+ \triangle DX_c \Lam_a+\triangle KX_c+ MX_c(\Lam_a^2-\Lam_c^2)+DX_c(\Lam_a-\Lam_c) \\
&=& MX_c [Z (X_c^\star MX_c \Lam_a+\epsilon_1 \epsilon_2 \Lam_c^\star X_c^\star MX_c+X_c^\star DX_c )\Lam_a+\Lam_c Z (X_c^\star MX_c \Lam_a+\epsilon_1 \epsilon_2 \Lam_c^\star X_c^\star MX_c \\ && +X_c^\star DX_c )+ \Lam_a^2-\Lam_c^2]+DX_c \left[Z \left(X_c^\star MX_c \Lam_a+\epsilon_1 \epsilon_2 \Lam_c^\star X_c^\star MX_c+X_c^\star DX_c \right)+  \Lam_a-\Lam_c  \right] \\
&=& MX_c \left[Z R \Lam_a+\Lam_c ZR+ \Lam_a^2-\Lam_c^2  \right]+DX_c \left[Z R+\Lam_a- \Lam_c \right]
\\ &=& 0,
\eeano
where the last equality follows by setting $Z= \left(\Lam_c-\Lam_a \right)R^{-1},$ $R=X_c^\star MX_c \Lam_a+\epsilon_1 \epsilon_2 \Lam_c^\star X_c^\star MX_c+X_c^\star DX_c$. Therefore $(X_c, \Lam_a)$ is an invariant pair of $Q_{\triangle}(\lam).$


Finally, note that $Q_{\triangle}(\lam)\in \Q_n(\star,\epsilon_1,\epsilon_2)$ if and only if $\triangle M=\epsilon_1  \triangle M^\star,\,\triangle D=\epsilon_2  \triangle D^\star,\,\triangle K=\epsilon_1  \triangle K^\star.$ Hence $Q_{\triangle}(\lam) \in \Q_n(\star,\epsilon_1,\epsilon_2)$ if $Z=\epsilon_1 Z^\star$ holds. Moreover, $R=S-X_c^\star MX_c (\Lam_c-\Lam_a)$ where $S=X_c^\star MX_c \Lam_c+\epsilon_1 \epsilon_2 \Lam_c^\star X_c^\star MX_c+X_c^\star DX_c$. By Corollary \ref{S_jProperty}, we have $S^\star =\epsilon_2 S$ and $S\Lam_c=\epsilon_1 \epsilon_2 \Lam_c^\star S,$ that is, $S \Lam_c=\epsilon_1 (S \Lam_c)^\star$. Then 
\beano 
&& Z= \epsilon_1 Z^\star \\
 \Leftrightarrow &&  R^\star (\Lam_c-\Lam_a)=\epsilon_1 (\Lam_c^\star -\Lam_a^\star )R \\
  \Leftrightarrow && \left[S^\star -\epsilon_1 (\Lam_c^\star -\Lam_a^\star )X_c^\star MX_c \right] (\Lam_c- \Lam_a)=\epsilon_1 (\Lam_c^\star -\Lam_a^\star ) \left[S- X_c^\star MX_c (\Lam_c-\Lam_a) \right]  \\
\Leftrightarrow && S^\star \Lam_c-S^\star \Lam_a=\epsilon_1 \Lam_c^\star S-\epsilon_1 \Lam_a^\star S \\
\Leftrightarrow && S\Lam_a=\epsilon_1(S\Lam_a)^\star,
\eeano where the last statement follows by using $S^\star =\epsilon_2 S$ and $S\Lam_c=\epsilon_1 \epsilon_2 \Lam_c^\star S$. Therefore the perturbed quadratic matrix polynomial $Q_{\triangle}(\lam) \in \Q_n(\star,\epsilon_1,\epsilon_2)$ whenever $S\Lam_a=\epsilon_1(S\Lam_a)^\star$ holds. Hence the desired result follows. $\hfill{\square}$
\end{proof}

\begin{remark}
Note that the spectral condition $(1)$ in Theorem \ref{P2_soln} is generically satisfied. Indeed, the spectrum of $\Lam_c$ is closed with respect to $(\star,\epsilon_1,\epsilon_2)$-symmetry, that is, if $\sigma(\Lam_c)= \sigma(\epsilon_1 \epsilon_2 \Lam_c^\star)$ then the condition $(1)$ in Theorem \ref{P2_soln} is satisfied if all the eigenvalues of $\Lam_c$ are different from the eigenvalues of $\Lam_f$.
\end{remark}

The next theorem describes a set of structured perturbations of a structured matrix polynomial that change the invariant pair $(X_c,\Lam_c)$ by $(X_a=X_cP,\Lam_a)$ while fixing $(X_f,\Lam_f)$ as an invariant pair of the perturbed polynomials, where $P$ is nonsingular. 

\begin{theorem} \label{X_a=X_cZ}
Suppose $(X_c,\Lam_c)\in \K^{n \times p_1} \times \K^{p_1 \times p_1}$ and $(X_f,\Lam_f)\in \K^{n \times p_2} \times \K^{p_2 \times p_2}$ are invariant pairs of a matrix polynomial $Q(\lam)=\lam^2 M+\lam D+ K \in  \Q_n(\st,\ep_1,\ep_2)$. Suppose $\Lam_a\in\K^{p_1\times p_1}$. Let $P\in \K^{p_1 \times p_1}$ be a nonsingular matrix such that $R:=X_c^\star MX_c P\Lam_a P^{-1}+\epsilon_1 \epsilon_2 \Lam_c^\star X_c^\star MX_c+X_c^\star DX_c$ is nonsingular. Let $\sig(\Lam_c)\cap \sigma(\epsilon_1 \epsilon_2 \Lam_f^\star) =\emptyset$. Then $(X_a=X_c P,\Lam_a)$ and $(X_f,\Lam_f)$ are invariant pairs of the quadratic matrix polynomial $Q_{\triangle}(\lam)=\lam^2 (M+\triangle M)+\lam (D+\triangle D)+(K+\triangle K)\in \K^{n \times n}[\lam]$ where
\beano 
\triangle M &=& MX_c Z X_c^\star M, \\
\triangle D &=& \epsilon_1 \epsilon_2 MX_c Z \Lam_c^\star X_c^\star M+MX_c Z X_c^\star D+MX_c \Lam_c Z X_c^\star M+DX_c Z X_c^\star M, \\
\triangle K &=& \epsilon_1 \epsilon_2 MX_c \Lam_c Z \Lam_c^\star X_c^\star M+MX_c \Lam_c Z X_c^\star D+\epsilon_1 \epsilon_2 DX_c Z \Lam_c^\star X_c^\star M+DX_c Z X_c^\star D,
\eeano 
$Z=(\Lam_c-P\Lam_a P^{-1})R^{-1}.$


Moreover, if  $SP\Lam_a P^{-1}=\epsilon_1 (SP \Lam_a P^{-1})^\star$ then $Q_{\triangle}(\lam) \in \Q_n(\star,\epsilon_1,\epsilon_2)$ where $S:=X_c^\star MX_c \Lam_c+\epsilon_1 \epsilon_2 \Lam_c^\star X_c^\star MX_c+X_c^\star DX_c$.

\end{theorem}
\begin{proof}
Since $P$ is nonsingular so it can be verified that $(X_a=X_cP,\Lam_a)$ is an invariant pair of $Q_{\triangle}(\lam)$ if and only if $(X_c,P\Lam_aP^{-1})$ is an invariant pair of $Q_{\triangle}(\lam)$. Hence the desired result follows from Theorem \ref{P2_soln}. $\hfill{\square}$
\end{proof}


Observe that the assumption of the existence of a nonsingular matrix $P\in \K^{p_1 \times p_1}$ in Theorem \ref{X_a=X_cZ} such that $R=X_c^\star MX_c P\Lam_a P^{-1}+\epsilon_1 \epsilon_2 \Lam_c^\star X_c^\star MX_c+X_c^\star DX_c$ is nonsingular, is not obvious. Indeed, there exists a nonsingular matrix $P$ such that $R$ will be nonsingular if and only if the matrix equation $$(R-R_1)P - (X_c^\st MX_c)P\Lam_a =0,$$ where $R_1=\epsilon_1\epsilon_2\Lam_c^\st (X_c^\st MX_c) + (X_c^\st DX_c)$ has a nonsingular solution $P$ for a nonsingular matrix $R.$ Setting $R=SP^{-1}$ for some nonsingular matrix $S,$ the matrix equation reduces to \begin{equation}\label{eqn:sylves1} R_1P + (X_c^\st MX_c)P\Lam_a = S,\end{equation} and the problem is to find a nonsingular solution $P$ of this equation for a given nonsingular matrix $S.$ If $X_c^\st MX_c$ is nonsingular then the above equation becomes the Sylvester equation of the form \begin{equation}\label{eqn:sylves}TA-BT=C,\end{equation} where $T=P,$ $B=-(X_c^\st MX_c)^{-1}R_1,$ $A=\Lam_a$, and $C=(X_c^\st MX_c)^{-1}S,$ which is nonsingular. Now the equation (\ref{eqn:sylves}) has a nonsingular solution $T$ if $\mathfrak{R}(C)\subseteq \mathfrak{R}(T)$ and $(B,C)$ is controllable, that is, $\bmatrix{C & BC & \hdots & B^{p_1-1}C}$ has rank $p_1$ [Theorem 2, \cite{hearon1977nonsingular}]. Now since $C$ is nonsingular, it is evident that these conditions are satisfied, and hence such a nonsingular matrix $P$ can be obtained that satisfies equation (\ref{eqn:sylves1}).

We emphasize that the expressions of the perturbation matrices $\triangle M,\,\triangle D,\,\triangle K$ in Theorem $\ref{P2_soln}$ and Theorem \ref{X_a=X_cZ} do not require knowledge of $(X_f,\Lam_f).$

\section{Solution of SEEP for quadratic matrix polynomials with symmetry structures}

Let $Q(\lam)=\lam^2 M+\lam D + K\in \Q_n(\st,\epsilon_1,\epsilon_2) \subseteq\K^{n\times n}[\lam]$ be a quadratic matrix polynomial. Then due to the symmetry structures of the coefficients, the polynomial $Q(\lam)$ inherits certain spectral symmetry, called eigenvalue-paring \cite{mackey2006structured}. Besides,  the algebraic, geometric, and partial multiplicities of the two eigenvalues in the pair are equal. Indeed, from Proposition \ref{Property} (a) it follows that if $\lam_0$ is an eigenvalue of $Q(\lam)$ then $\epsilon_1\epsilon_2\lam^\st_0$ is also an eigenvalue of $Q(\lam)$, when $\lam_0\neq \epsilon_1\epsilon_2\lam^\st_0.$ This induces an eigenvalue pairing $(\lam_0,\epsilon_1\epsilon_2\lam^\st_0)$ of $Q(\lam).$ Then the SEEP for $(\st,\epsilon_1,\epsilon_2)$-structured polynomials is described as follows \cite{Chu_postdefn,Datta,Mao_Dai}. Let $\sigma^c_s=\{\lam_j^c,\, \epsilon_1\epsilon_2(\lam_j^c)^\st,\, \lam_k^c : \lam_j^c\neq \epsilon_1\epsilon_2(\lam_j^c)^\st,\, \lam_k^c=\epsilon_1\epsilon_2(\lam_k^c)^\st,\, j=1,\hdots,m,\, k=2m+1,\hdots,p\}$ be a set of finite eigenvalues of $Q(\lam)$ that are known, and $\lam_l^f,\, l=p+1, \hdots, 2n$ are the rest of the eigenvalues that need not be known. Then, given a set of scalars $\sigma^a_s=\{\lam_j^a,\, \epsilon_1\epsilon_2(\lam_j^a)^\st,\, \lam_k^a : \lam_j^a \neq \epsilon_1\epsilon_2(\lam_j^a)^\st,\, \lam_k^a=\epsilon_1\epsilon_2(\lam_k^a)^\st, \,j=1,\hdots,m,\, k=2m+1,\hdots,p\},$ determine  polynomials $\triangle Q(\lam)=\lam^2 \triangle M+\lam \triangle D + \triangle K\in \Q_n(\st,\epsilon_1,\epsilon_2)$ such that the complete set of eigenvalues of the perturbed polynomial $Q_\triangle (\lam)=\lam^2 (M+\triangle M) + \lam (D+\triangle D) + (K+\triangle K)$ is given by $\{\lam^a_j,\, \epsilon_1\epsilon_2(\lam_j^a)^\st,\, \lam_k^a,\, \lam^f_l : 1\leq j\leq m,\, 2m+1\leq k\leq p,\, p+1\leq l\leq 2n\}.$ 

Then the SEEP can be defined in terms of preserving invariant pairs of structured polynomials under structure-preserving perturbations. The invariant pairs are defined by using eigenpairs of the polynomial. Let $(X_c,\Lam_c)$ and $(X_f,\Lam_f)$ be the invariant pairs of a given polynomial $Q(\lam)\in\Q_n(\st,\epsilon_1,\epsilon_2)\subset \K^{n\times n}[\lam]$, where $\Lam_c,\,\Lam_f$ are block diagonal matrices whose eigenvalues are eigenvalues of $Q(\lam)$ and the columns of $X_c,\, X_f$ are eigenvectors corresponding to those eigenvalues respectively. Then given a block diagonal matrix $\Lam_a$ determine polynomials $\triangle Q(\lam)=\lam^2 \triangle M + \lam \triangle D + \triangle K\in\K^{n\times n}[\lam]$ such that $$Q_\triangle(X_cP,\Lam_a)=0 \,\, \mbox{and} \,\, Q_\triangle (X_f,\Lam_f)=0,$$ where $Q_\triangle(\lam)=Q(\lam)+\triangle Q(\lam)\in \Q_n(\st,\epsilon_1,\epsilon_2)$ and $P$ is a suitably chosen nonsingular matrix.  Note that the invariant pair $(X_f,\Lam_f)$ need not be known, and $\sigma(\Lam_c)\cap \sigma(\Lam_f)=\emptyset.$ Below we describe the structure of the invariant pairs of $Q(\lam)$ that depends on the structure of the polynomial $Q(\lam).$
\begin{itemize}
    \item \textbf{$Q(\lam)\in \Q_n(*,\epsilon_1,\epsilon_2)\subset \C^{n\times n}[\lam]:$} Let $x^c_j,\,\widetilde{x}^c_j,$ and $x_k^c$ denote the eigenvectors corresponding to the eigenvalues $\lam_j^c,$ $\epsilon_1\epsilon_2 \overline{\lam_j^c}\, (\neq \lam_j^c)$ and $\lam_k^c=\epsilon_1\epsilon_2 \overline{\lam_k^c}$ of $Q(\lam)$ respectively, $j=1,\hdots,m,\, k=2m+1,\hdots,p$. Let $\lam^a_j,\,\ep_1 \ep_2 \overline{\lam^a_j},\,\lam^a_k,\,j=1,\hdots,m,\,k=2m+1,\hdots,p$ be a collection of scalars such that $\lam^a_j \neq \ep_1 \ep_2 \overline{\lam^a_j}$ and $\lam^a_k=\ep_1 \ep_2 \overline{\lam^a_k}$. Further, let $x_l^f$ denote an eigenvector corresponding to the eigenvalue $\lam_l^f,$ $l=p+1,\hdots, 2n$ which need not be known. Then setting $$\Lam_j^c=\bmatrix{\lam_j^c & 0\\ 0 & \epsilon_1\epsilon_2 \overline{\lam^c_j}}, \Lam^a_j=\bmatrix{\lam_j^a & 0\\ 0 & \epsilon_1\epsilon_2 \overline{\lam^a_j}}, X_j^c=\bmatrix{x_j^c & \widetilde{x}_j^c},\, j=1,\hdots,m,$$ we have $$Q(X_c,\Lam_c) =0 \,\, \mbox{where} 
\left\{
  \begin{array}{l}
\Lam_c=\diag (\Lam_1^c,\hdots, \Lam_m^c,\lam_{2m+1}^c, \hdots, \lam_p^c),  \\
X_c=\bmatrix{X_1^c & \hdots & X_m^c & x_{2m+1}^c & \hdots & x_p^c}. 
  \end{array}
\right.
$$ Also, $\Lam_f=\diag(\lam_{p+1}^f, \hdots, \lam_{2n}^f),$ $X_f=\bmatrix{x_{p+1}^f & \hdots & x_{2n}^f}$ such that $Q(X_f,\Lam_f)=0,$ and $\Lam_a=\diag (\Lam_1^a,\hdots, \Lam_m^a,\lam_{2m+1}^a, \hdots, \lam_p^a).$

\item $Q(\lam)\in \Q_n(T,1,1)\subset \C^{n\times n}[\lam]:$ Let $(\lam^c_j,x_j^c),\, 1\leq j\leq p$ and $(\lam_l^f,x_l^f),\, p+1\leq l\leq 2n$ be the collection of eigenpairs of $Q(\lam)$, where the later pairs need not be known. Then we have $Q(X_c,\Lam_c)=0$ where $$\Lam_c=\diag(\lam_1^c,\hdots,\lam_p^c), X_c=\bmatrix{x_1^c &\hdots & x_p^c}$$ and $\Lam_a=\diag(\lam^a_1,\hdots,\lam_p^a)$ where $\lam_j^a$s are scalars, $\Lam_f=\diag(\lam_{p+1}^f, \hdots, \lam_{2n}^f),$ $X_f=\bmatrix{x_{p+1}^f & \hdots & x_{2n}^f}$ with $Q(X_f,\Lam_f)=0.$

 \item \textbf{$Q(\lam)\in \Q_n(T,1,1)\subset \R^{n\times n}[\lam]:$} Let $(\lam_j^c, x_j^c),$ $(\overline{\lam_j^c}, \overline{x_j^c})$ and $(\lam_k^c, x_k^c)$ are known eigenpairs of $Q(\lam),$ where $\lam_j^c\in\C \smallsetminus\R,$ $\lam^c_k \in\R$, $1\leq j\leq m,$ $2m+1\leq k\leq p.$ Let $\lam_l^f,$ $p+1\leq l\leq 2n$ be the rest of the eigenvalues of $Q(\lam)$ and $x^f_l$ are the corresponding eigenvectors, which need not be known. Let $\lam^a_j,\,\lam^a_k,\,j=1,\hdots,m,\,k=2m+1,\hdots,p$ be a collection of scalars such that $\lam^a_j\in \C \smallsetminus \R,\,\lam^a_k \in \R$. Then setting $$\Lam_j^c=\bmatrix{\re(\lam_j^c) & \im(\lam_j^c)\\ -\im(\lam_j^c) & \re(\lam_j^c)},\, \Lam^a_j=\bmatrix{\re(\lam_j^a) & \im(\lam_j^a)\\ -\im(\lam_j^a) & \re(\lam_j^a)},\, X_j^c=\bmatrix{\re(x_j^c) & \im(x_j^c)},\, j=1,\hdots,m,$$ we have $$Q(X_c,\Lam_c) =0 \,\, \mbox{where} 
\left\{
  \begin{array}{l}
\Lam_c=\diag (\Lam_1^c,\hdots, \Lam_m^c,\lam_{2m+1}^c, \hdots, \lam_p^c),  \\
X_c=\bmatrix{X_1^c & \hdots & X_m^c & x_{2m+1}^c & \hdots & x_p^c}. 
  \end{array}
\right.
$$ Also, $\Lam_a=\diag (\Lam_1^a,\hdots, \Lam_m^a,\lam_{2m+1}^a, \hdots, \lam_p^a).$ Besides, $\Lam_f=\diag(\lam_{p+1}^f, \hdots, \lam_{2n}^f),$ \\ $X_f=\bmatrix{x_{p+1}^f & \hdots & x_{2n}^f}$ such that $Q(X_f,\Lam_f)=0.$

 \item \textbf{$Q(\lam)\in \Q_n(T,\epsilon_1,-\epsilon_1)\subset \C^{n\times n}[\lam]:$} Let $(\lam_j^c, x_j^c),$ $(-\lam_j^c,\widetilde{x}_j^c),\,j=1,\hdots,p$ be known eigenpairs of $Q(\lam),$ where $0\neq \lam_j^c\in\C.$ Let $\lam_l^f,\, 2p+1\leq l\leq 2n$ be the rest of the eigenvalues of $Q(\lam),$ and $x_l^f$ are corresponding eigenvectors, which need not be known. Let $\lam^a_j,\,j=1,\hdots,p$ be a collection of scalars. Then setting $$\Lam_j^c=\bmatrix{\lam_j^c & 0\\ 0 & -\lam^c_j},\, \Lam^a_j=\bmatrix{\lam_j^a & 0\\ 0 & -\lam^a_j},\, X_j^c=\bmatrix{x_j^c & \widetilde{x}_j^c},\, j=1,\hdots,p,$$ we have $$Q(X_c,\Lam_c) =0 \,\, \mbox{where} 
\left\{
  \begin{array}{l}
\Lam_c=\diag (\Lam_1^c,\hdots, \Lam_p^c),  \\
X_c=\bmatrix{X_1^c & \hdots & X_p^c}. 
  \end{array}
\right.
$$ Also, $\Lam_f=\diag(\lam_{2p+1}^f, \hdots, \lam_{2n}^f),$ $X_f=\bmatrix{x_{2p+1}^f & \hdots & x_{2n}^f}$ such that $Q(X_f,\Lam_f)=0,$ and $\Lam_a=\diag (\Lam_1^a,\hdots, \Lam_p^a).$

 \item \textbf{$Q(\lam)\in \Q_n(T,\epsilon_1,-\epsilon_1)\subset \R^{n\times n}[\lam]:$} Let $(\lam^c_j,x^c_j),\,(\overline{\lam^c_j},\overline{x^c_j}),\, (-\lam^c_j,\tilde{x}^c_j),\,(-\overline{\lam^c_j},\overline{\tilde{x}^c_j}),$ $(\lam^c_k,x^c_k),\,(\overline{\lam^c_k},\overline{x^c_k}),$ $(\lam^c_l,x^c_l),\,(-\lam^c_l,\tilde{x}^c_l)$ be eigenpairs of $Q(\lam)$ with $\lam^c_j \neq 0 \in \C \smallsetminus (\R \cup \mathrm{i}\R),$ $\lam^c_k \neq 0 \in \mathrm{i}\R,$ $\lam^c_l \neq 0 \in \R,$ $j=1,\hdots,m_1,\,k=m_1+1,\hdots,m_2,\,l=m_2+1,\hdots,p$. Let $\lam^f_r,\,r=2m_1+2p+1,\hdots, 2n$ be the rest of the eigenvalues of $Q(\lam)$ corresponding to eigenvectors $x_r^f$, which need not be known. Let $\lam^a_j,\,\lam^a_k,\,\lam^a_l,\,j=1,\hdots,m_1,\,k=m_1+1,\hdots,m_2,\,l=m_2+1,\hdots,p$ be a collection of scalars where $\lam^a_j \in \C \smallsetminus (\R \cup \mathrm{i}\R),\,\lam^a_k \in \mathrm{i}\R$ and $\lam^a_l  \in \R$. Then setting 
 \beano 
&& \tilde{\Lam}^c_j=\bmatrix{\re(\lam^c_j)& \im(\lam^c_j) \\ -\im(\lam^c_j)& \re(\lam^c_j)},\,\Lam^c_k=\bmatrix{0& \im(\lam^c_k) \\ -\im(\lam^c_k)& 0},\,\Lam^c_l=\bmatrix{\lam^c_l & 0 \\ 0 & -\lam^c_l}, \\ &&\tilde{\Lam}^a_j=\bmatrix{\re(\lam^a_j)& \im(\lam^a_j) \\ -\im(\lam^a_j)& \re(\lam^a_j)},\,\Lam^a_k=\bmatrix{0& \im(\lam^a_k) \\ -\im(\lam^a_k)& 0},\,\Lam^a_l=\bmatrix{\lam^a_l & 0 \\ 0 & -\lam^a_l},
 \eeano we have 
 $$Q(X_c,\Lam_c) =0 \,\, \mbox{where} 
\left\{
  \begin{array}{l}
\Lam_c=\diag (\Lam_1^c,\hdots, \Lam_{m_1}^c, \Lam^c_{m_1+1}, \hdots, \Lam_{m_2}^c, \Lam_{m_2 + 1}^c, \hdots, \Lam_p^c), \\
X_c=\bmatrix{X_1^c & \hdots & X_{m_1}^c\,\, X_{m_1+1}^c & \hdots & X_{m_2}^c\,\, X_{m_2+1}^c & \hdots & X_p^c},
  \end{array}
\right.
$$ where $\Lam^c_j=\diag(\tilde{\Lam}^c_j,\, -\tilde{\Lam}^c_j)$, $X^c_j=\bmatrix{G^c_j& H^c_j},\,G^c_j=\bmatrix{\re(x^c_j)& \im(x^c_j)},\,H^c_j=\bmatrix{\re(\tilde{x}^c_j)& \im(\tilde{x}^c_j)},\\ X^c_k=\bmatrix{\re(x^c_k)& \im(x^c_k)},$ $X^c_l=\bmatrix{x^c_l& \tilde{x}^c_l}.$ Also, $$\Lam_a=\mathrm{diag}\left(\Lam^a_1,\hdots,\Lam^a_{m_1},\,\Lam^a_{m_1+1},\hdots, \Lam^a_{m_2},\,\Lam^a_{m_2+1},\hdots, \Lam^a_p \right),$$ where $\Lam^a_j=\diag(\tilde{\Lam}^a_j,\, -\tilde{\Lam}^a_j),\,j=1,\hdots,m_1,$ and $\Lam_f=\diag(\lam_{2m_1+2p+1}^f, \hdots, \lam_{2n}^f),$ $\linebreak X_f=\bmatrix{x_{2m_1 + 2p+1}^f & \hdots & x_{2n}^f}$ with $Q(X_f,\Lam_f)=0.$

\end{itemize}

Then analytical solutions of MUP can be obtained by employing Theorem \ref{thm:mupip} for structured polynomials. The following corollary provides analytical solutions of SEEP for $(\st,\epsilon_1,\epsilon_2)$-structured polynomials by employing Theorem \ref{X_a=X_cZ}.  


\begin{corollary}\label{Cor:main}
Let $(X_c,\Lam_c)$ and $(X_f,\Lam_f)$ be the invariant pairs corresponding to eigenpairs of a polynomial $Q(\lam)\in\Q_n(\st, \epsilon_1,\epsilon_2)$ as discussed above and the later pair need not be known. Let the eigenvalues of $\Lam_c$ be simple and distinct. If $\sigma(\Lam_c)\cap \sigma(\epsilon_1\epsilon_2\Lam_f^\st)=\emptyset$ then the polynomial $\triangle Q(\lam) = \lam^2 \triangle M + \lam \triangle D +\triangle K\in \Q_n(\st, \epsilon_1, \epsilon_2)$ such that $Q_\triangle(X_a=X_cP,\Lam_a)=0=Q_\triangle(X_f,\Lam_f)$ where $Q_\triangle(\lam)=Q(\lam)+\triangle Q(\lam)$ is given by Theorem \ref{X_a=X_cZ}, and $P$ is an invertible matrix that depends on $\st\in\{*,T\},$ $\epsilon_1,\epsilon_2\in\{1, -1\}$ such that $R=X_c^\star MX_c P\Lam_a P^{-1}+\epsilon_1 \epsilon_2 \Lam_c^\star X_c^\star MX_c+X_c^\star DX_c$ is nonsingular. The solution matrix $P$ which defines the structured polynomial $\triangle Q(\lam)$ is constructed as:
\begin{itemize}
    \item \textbf{$Q(\lam)\in \Q_n(*,\epsilon_1,\epsilon_2)\subset \C^{n\times n}[\lam]:$} $P=\diag(P_1, \hdots, P_m, I_{p-2m})$ where  $$P_j=\bmatrix{\sqrt{\ep_1} \alpha_j a_j& (\lam^a_j-\ep_1 \ep_2 \overline{\lam^a_j})^{-1} b_j \\ -\ep_1(\lam^a_j-\ep_1 \ep_2 \overline{\lam^a_j})^{-1} a_j& \sqrt{\ep_1}  b_j \overline{\alpha}_j},$$ $\alpha_j=2\ep_1 \ep_2 \overline{\lam^c_j} (x^c_j)^*M\tilde{x}^c_j+(x^c_j)^*D\tilde{x}^c_j,$ and $a_j, b_j, 1\leq j\leq m$ are arbitrarily chosen nonzero real numbers. 
    
      \item \textbf{$Q(\lam)\in \Q_n(T, 1, 1)\subset \C^{n\times n}[\lam]:$} $P$ is a diagonal matrix of order $p\times p$ with non-zero diagonal entries.
      
       \item \textbf{$Q(\lam)\in \Q_n(T, 1, 1)\subset \R^{n\times n}[\lam]:$} $P=\diag(P_1, \hdots, P_m, I_{p-2m}),$ where
       $$P_j=\left\{
  \begin{array}{l}
\bmatrix{a_j& \dfrac{\beta_j a_j}{\alpha_j} - \dfrac{\alpha_j}{2} \\ \dfrac{\beta_j a_j}{\alpha_j} + \dfrac{\alpha_j}{2}& -a_j}, \,\, \mbox{if} \,\, \alpha_j\neq 0 \\
\bmatrix{a_j& 1\\ -a_j& a_j^2}, \,\, \mbox{otherwise} 
  \end{array}
\right.$$ $\alpha_j=\re(\gamma_j)/2,\,\beta_j=-\im(\gamma_j)/2,$  $\gamma_j=2 \overline{\lam^c_j} (x^c_j)^*M \overline{x^c_j}+(x^c_j)^*D \overline{x^c_j},$ and  $a_j$ is an arbitrarily chosen real number such that $P_j$ is nonsingular. 
  
 \item \textbf{$Q(\lam)\in \Q_n(T, \epsilon_1, -\epsilon_1)\subset \C^{n\times n}[\lam]:$} $$P=\left\{
  \begin{array}{l}
\diag(P_1, \hdots, P_p), \,\, \mbox{if} \,\, \epsilon_1=1 \\
W\otimes I_2, \,\, \mbox{if} \,\, \epsilon_1=-1, 
  \end{array}
\right.$$ where $P_j\in\C^{2\times 2}$ is an arbitrary nonsingular matrix, and $W$ is a diagonal matrix of order $p\times p$ with non-zero diagonal entries, $\otimes$ denotes the Kronecker product of matrices.     
 
 \item \textbf{$Q(\lam)\in \Q_n(T, \epsilon_1, -\epsilon_1)\subset \R^{n\times n}[\lam]:$} Set $\gamma_j=-2\overline{\lam^c_j} (x^c_j)^*M \overline{\tilde{x}^c_j}+(x^c_j)^*D \overline{\tilde{x}^c_j},\,\alpha_j=\re(\gamma_j)/2,\,\beta_j=-\im(\gamma_j)/2$. Then $$P=\left\{
  \begin{array}{l}
\diag(P_1, \hdots, P_{m_1}, P_{m_1+1}, \hdots, P_{m_2}, P_{m_2+1}, \hdots, P_{p}), \,\, \mbox{if} \,\, \epsilon_1=1, \\
\diag(P_1, \hdots, P_{m_1}, I_{2p-2m_1}), \,\, \mbox{if} \,\, \epsilon_1=-1, 
  \end{array}
\right.$$ where 
\beano 
P_j &=& 
\left\{
  \begin{array}{l}
  \bmatrix{0& \mathfrak{P}_j\\-\mathfrak{P}_j& 0}, \,\, \mbox{if} \,\, \epsilon_1=1, \\
  \bmatrix{0& \mathfrak{P}_j\\ \mathfrak{P}_j& 0}, \,\, \mbox{if} \,\, \epsilon_1=-1,
  \end{array}
\right.
 \mathfrak{P}_j=
\left\{
  \begin{array}{l}
\bmatrix{r_j& \frac{\beta_jr_j}{\alpha_j}-\frac{\alpha_j}{2} \\ \frac{\beta_jr_j}{\alpha_j}+\frac{\alpha_j}{2}& -r_j} \,\, \mbox{if} \,\, \alpha_j\neq 0, \\
\bmatrix{r_j& 1\\-r_j& r_j^2}, \,\, \mbox{if} \,\, \alpha_j= 0, 
  \end{array}
\right. r_j\in \R,\, 1\leq j\leq m_1 \\
P_k &=& \bmatrix{a_k & b_k \\ -c_k & c_k},\,\, a_k,\,b_k,\, c_k\in \R,\, m_1+1\leq k \leq m_2, 
\eeano and $P_l \in \R^{2\times 2},\, m_2+1\leq l\leq p.$ The free parameters are chosen such that $P_j,\, P_k,\, P_l$ are nonsingular matrices for all $j,\, k,\, l.$
\end{itemize}
\end{corollary}

\pf The proof follows from the fact that $SP\Lam_aP^{-1}=\epsilon_1 (SP\Lam_a P^{-1})^\st$ due to Theorem \ref{X_a=X_cZ} and Proposition \ref{Property}, where $S=X_c^\star MX_c \Lam_c+\epsilon_1 \epsilon_2 \Lam_c^\star X_c^\star MX_c+X_c^\star DX_c$. Indeed, $MX_c \Lam_c^2+DX_c \Lam_c+KX_c=0$ holds. Since $\sigma(\Lam_c)\cap \sigma(\epsilon_1\epsilon_2\Lam_f^\st)=\emptyset$ and nonsingular matrix $P$ is chosen such that $R$ is invertible then the expression for structured perturbation matrices $\triangle M,\,\triangle D,\,\triangle K$ follows from Theorem \ref{X_a=X_cZ}.

Let $Q(\lam)\in\Q_n(*,\epsilon_1,\epsilon_2).$ 
As the eigenvalues $\lam^c_j,\,\ep_1 \ep_2 \overline{\lam^c_j},\,\lam^c_k$ are distinct, applying Proposition \ref{Property} $(c)$ we obtain $S=X_c^*MX_c \Lam_c+\ep_1 \ep_2 \Lam_c^* X_c^*MX_c+X_c^*DX_c=\diag(S_1,\hdots, S_m,\,s_{2m+1},\hdots, s_p)$ where $S_j=(X^c_j)^*MX^c_j \Lam^c_j+\ep_1 \ep_2 (\Lam^c_j)^* (X^c_j)^*MX^c_j+(X^c_j)^*DX^c_j,\,s_k=2 \lam^c_k (x^c_k)^*Mx^c_k+(x^c_k)^*Dx^c_k,\,j=1,\hdots,m,\,k=2m+1,\hdots ,p$. Again by Corollary \ref{blockform} it follows that $$S_j=(X^c_j)^*MX^c_j \Lam^c_j+\ep_1 \ep_2 (\Lam^c_j)^* (X^c_j)^*MX^c_j+(X^c_j)^*DX^c_j=\bmatrix{0& \alpha_j\\ \ep_2 \overline{\alpha}_j& 0}$$ where $\alpha_j= 2\ep_1 \ep_2 \overline{\lam^c_j} (x^c_j)^*M \tilde{x}^c_j+(x^c_j)^*D \tilde{x}^c_j$. Then choosing $P=\mbox{diag}\left(P_1,\hdots, P_m,\,I_{p-2m} \right)$ with $$P_j=\bmatrix{\sqrt{\ep_1} \alpha_j a_j& (\lam^a-\ep_1 \ep_2 \overline{\lam^a_j})^{-1} b_j \\ -\ep_1(\lam^a-\ep_1 \ep_2 \overline{\lam^a_j})^{-1} a_j& \sqrt{\ep_1}  b_j \overline{\alpha}_j}$$ it follows that $SP\Lam_a P^{-1}=\ep_1(SP\Lam_a P^{-1})^*$ holds, where $a_j,\,b_j$ are arbitrarily chosen nonzero real numbers. Hence, by Theorem \ref{X_a=X_cZ} it follows that $Z=\ep_1 Z^*$ holds, thus $\triangle M =\ep_1 \triangle M^*,\,\triangle D =\ep_2 \triangle D^*,\,\triangle K =\ep_1 \triangle K^*$ holds. Therefore the updated quadratic matrix polynomial $Q_{\triangle}(\lam)=\lam^2(M+\triangle M)+\lam(D +\triangle D)+(K+\triangle K) \in \Q_n(*,\ep_1,\ep_2)$. The results for other structures follow using similar arguments.
$\hfill{\square}$

Next we present the following remark. 
\begin{remark} \label{Reco_Chu_postdefn} (Recovery of results in Chu et al. \cite{Chu_postdefn}, and  Kuo and Datta \cite{Kuo_Datta})
Let $Q(\lam)=\lam^2 M + \lam D + K$ be a polynomial such that $M, K$ are real symmetric positive definite matrices and $D$ is a real symmetric matrix. Thus the eigenvalues that are to be changed, are nonzero, and hence $\Lam_c$ is nonsingular. Now $Q(X_c,\Lam_c)=0$ implies $MX_c \Lam_c+DX_c=-KX_c \Lam_c^{-1}$ and $\Lam_c^TX_c^TM+X_c^TD=-(\Lam_c^T)^{-1} X_c^T K$. Assume that eigenvalues of $\Lam_a$ are nonzero. Thus by Corollary \ref{Cor:main} when $Q(\lam)\in \Q_n(T,1,1)\subset\R^{n\times n}[\lam]$ we have
\beano
\triangle M &=& MX_c Z X_c^TM, \\ 
\triangle D &=& MX_c Z \left(\Lam_c^TX_c^TM+X_c^TD \right)+\left(MX_c \Lam_c+DX_c\right) Z X_c^TM \\
&=& -MX_cZ (\Lam_c^T)^{-1}X_c^TK-KX_c \Lam_c^{-1} Z X_c^TM, \\ 
\triangle K &=& \left(MX_c \Lam_c+DX_c\right) Z \left(\Lam_c^TX_c^TM+X_c^TD \right)= KX_c\Lam_c^{-1} Z (\Lam_c^T)^{-1} X_c^TK,\eeano 
where \beano Z &=& \left(\Lam_c-P \Lam_aP^{-1} \right) \left(X_c^TMX_cP \Lam_a P^{-1}+(\Lam_c^TX_c^TM+X_c^TD)X_c \right)^{-1} \\
&=& \left(\Lam_c-P \Lam_aP^{-1} \right) \left(X_c^TMX_cP \Lam_a P^{-1}-(\Lam_c^T)^{-1} X_c^TKX_c \right)^{-1} \\
&=& \left(P\Lam_a-\Lam_cP \right) \left( (\Lam_c^T)^{-1}X_c^TKX_cP-X_c^TMX_cP \Lam_a \right)^{-1} \\
&=& \left(P\Lam_a-\Lam_cP \right) \left(X_c^TKX_cP-\Lam_c^T X_c^TMX_cP \Lam_a \right)^{-1} \Lam_c^T.
\eeano  

Besides, choosing $P$ as described in Corollary \ref{Cor:main} in such a way that $((P\Lam_aP^{-1})^TX_c^TMX_c-X_c^TKX_c \Lam_c^{-1} )(\Lam_c-P\Lam_aP^{-1})$ (which is symmetric due to the definition of $P$) is a positive semi-definite matrix, we obtain a positive semi-definite matrix $Z$ and consequently $\triangle M,$ $\triangle K$ are also symmetric positive semi-definite matrices. Thus $M+\triangle M$ and $K+\triangle K$ are real symmetric positive definite matrices. It can be verified that the structured perturbations proposed in \cite{Chu_postdefn} and the obtained perturbations are same.


On the other hand, $M$ is a nonsingular symmetric matrix, and $D,\,K$  are symmetric matrices in \cite{Kuo_Datta}. The perturbations proposed in \cite{Kuo_Datta} are given by
\beano
\triangle M &=& (d-1)M+MX_c \Phi X_c^T M, \\
\triangle D &=& (d-1)D-MX_c \Phi (\Lam_c^T)^{-1}X_c^TK-KX_c \Lam_c^{-1} \Phi X_c^TM, \\
\triangle K &=& (d-1)K+KX_c \Lam_c^{-1} \Phi (\Lam_c^T)^{-1} X_c^TK
\eeano
where $d \neq 0 \in \R$ and $\Phi$ is a symmetric matrix that can be obtained by solving the Sylvester equation 
$$\left(\Theta^T M_1-K_1\Lam_c^{-1}\right) \Phi M_1+M_1 \Phi \left(M_1 \Theta-(\Lam_c^T)^{-1}K_1\right)=d\left(\Lam_c-\Theta \right)^TM_1+d M_1\left(\Lam_c-\Theta \right)$$
with $M_1=X_c^TMX_c,\,K_1=X_c^TKX_c,\,\Theta=P \Lam_a P^{-1}$. Note that these expressions coincide with the expressions obtained in this paper as given above by setting $d=1,$ and $\Phi=Z.$  

\end{remark}

It may further be noted that the matrix $\Lam_c$ is considered as a nonsingular matrix in \cite{Kuo_Datta}, that is, the proposed perturbations can change only the nonzero eigenvalues. 

Now summarizing the above results we provide an algorithm for solving SEEP that arises in structural models, that is, for quadratic polynomials $Q(\lam)=\lam^2 M+\lam D+K\in\R^{n\times n}[\lam]$ with $M, K$ are positive definite matrices, and $D$ is a symmetric matrix.


\noindent\rule{14.2cm}{0.5pt} \\
\textbf{Algorithm } \\
\noindent\rule{14.2cm}{0.5pt} \\
\textbf{Input:} Real symmetric positive definite matrices $M,\,K$ and real symmetric matrix $D$. \\
\textbf{Output:} Real symmetric positive semi-definite matrices $\triangle M,\,\triangle K$ and real symmetric matrix $\triangle D$. \\
\textbf{1.} Form the matrices $X_c,\,\Lam_c,\,\Lam_a$ as mentioned above. \\ 
\textbf{2.} Choose the matrix $P$ as defined in Corollary \ref{Cor:main} for the case when $Q(\lam)=\lam^2 M+\lam D+K \in \Q_n(T,1,1) \subset \R^{n \times n}[\lam]$ in such a way that $((P\Lam_aP^{-1})^TX_c^TMX_c-X_c^TKX_c \Lam_c^{-1} )(\Lam_c-P\Lam_aP^{-1})$ is positive semi-definite.\\
\textbf{3.} Compute the matrix $Z$ as defined in Remark \ref{Reco_Chu_postdefn}. Consequently we obtain the positive semi-definite matrices $\triangle M,\,\triangle K$ and symmetric matrix $\triangle D$ by applying Remark \ref{Reco_Chu_postdefn}.  \\
\noindent\rule{14.2cm}{0.5pt} \\

Next we have the following remark.
\begin{remark} \label{Reco_Mao_Dai} (Recovery of results in Mao and Dai \cite{Mao_Dai})
Let $Q(\lam)=\lam^2 M + \lam D+K\in\R^{n\times n}[\lam]$ be a matrix polynomial with $M,\,K$ as symmetric positive definite matrices and $D$ is a skew-symmetric matrix. Then $Q(\lam)$ has purely imaginary eigenvalues. Thus $(\lam_0,x_0)$ is an eigenpair of $Q(\lam)$ if and only if $(\overline{\lam}_0,\overline{x}_0)$ is also an eigenpair of $Q(\lam)$ where $0\neq \lam_0  \in \mathrm{i}\R$. Then the matrix $\Lam_c=\diag(\Lam^c_1,\hdots, \Lam^c_p)$ is nonsingular, where $\Lam^c_j=\bmatrix{0& \im(\lam^c_j)\\ -\im(\lam^c_j)& 0},\,j=1,\hdots,p.$ Since $\Lam_c$ is invertible so we have $MX_c \Lam_c+DX_c=-KX_c \Lam_c^{-1}$ and $\Lam_c^TX_c^TM-X_c^TD=-(\Lam_c^T)^{-1} X_c^TK$. Let $\lam^a_j \neq 0 \in \mathrm{i}\R,\,j=1,\hdots,p.$ Then set $\Lam_a=\diag(\Lam^a_1,\hdots, \Lam^a_p)$ with $\Lam^a_j=\bmatrix{0& \im(\lam^a_j)\\ -\im(\lam^a_j)& 0}$. 

Therefore by Corollary \ref{Cor:main} when $Q(\lam) \in \Q_n(T,1,-1) \subset \R^{n \times n}[\lam]$, we have
\beano 
\triangle M &=& MX_c Z X_c^TM, \\ 
\triangle D &=& -MX_c Z \left(\Lam_c^TX_c^TM-X_c^TD \right)+\left(MX_c \Lam_c+DX_c \right) Z X_c^TM \\
&=& MX_c Z (\Lam_c^T)^{-1} X_c^TK-KX_c \Lam_c^{-1} Z X_c^TM, \\ 
\triangle K &=& -\left(MX_c \Lam_c+DX_c \right) Z \left(\Lam_c^TX_c^TM-X_c^TD \right) = -KX_c\Lam_c^{-1} Z (\Lam_c^T)^{-1} X_c^TK \,\, 
\eeano with \beano Z &=& \left(\Lam_c-P \Lam_a P^{-1} \right) \left(X_c^TMX_cP \Lam_a P^{-1}-(\Lam_c^TX_c^TM-X_c^TD)X_c \right)^{-1} \\ &=& \left(\Lam_c P-P \Lam_a \right) \left(X_c^TMX_c P \Lam_a+(\Lam_c^T)^{-1} X_c^T KX_c P \right)^{-1}\eeano
where $P=\diag(P_1,\hdots,P_p),$ $P_j=\bmatrix{a_j& b_j \\-c_j& c_j},\,j=1,\hdots,p$ and $a_j,\,b_j,\,c_j$ are arbitrarily chosen real numbers for which $P$ and $R$ in Corollary \ref{Cor:main} are nonsingular matrices. 
Note that for such choice of $P,$ the matrix $Z$ is real symmetric and $$Z=\left(P^TX_c^TKX_c \Lam_c^{-1}-\Lam_a P^TX_c^TMX_c \right)^{-1} \left(\Lam_a P^T-P^T \Lam_c \right)$$ since $\Lam_c^T=-\Lam_c,\,\Lam_a^T=-\Lam_a.$ 

Observe that this solution realizes the solution obtained by Mao and Dai in \cite{Mao_Dai} as follows. The perturbations obtained in their paper (Theorem 3.1 and Theorem 3.2, \cite{Mao_Dai}) are given by
\beano
\triangle M &=& MX_cEX_c^TM, \\
\triangle D &=& MX_c E (\Lam_c^T)^{-1} X_c^TK-KX_c \Lam_c^{-1} E X_c^TM,\\
\triangle K &=& -KX_c\Lam_c^{-1} E (\Lam_c^T)^{-1} X_c^TK
\eeano
with $$E=\left(P^T X_c^TKX_c \Lam_c^{-1}-\Lam_aP^TX_c^T MX_c \right)^{-1} \left(\Lam_aP^T-P^T \Lam_c \right).$$ 
\end{remark}

\section{Numerical examples}
In this section, we consider numerical examples of structured matrix polynomials and illustrate the applications of the obtained solutions for SEEP. Let $(X_c,\Lam_c)$ be an invariant pair of a structured polynomial $Q(\lam)=\lam^2 M +\lam D+K.$ Let $(X_f,\Lam_f)$ and $(X_cP,\Lam_a)$ be the invariant pairs of the updated matrix polynomial $Q_\triangle (\lam)=\lam^2(M+\triangle M)+\lam (D+ \triangle D)+(K+\triangle K)$. Then we define the relative residuals of $(X_f,\Lam_f)$ and $(X_cP,\Lam_a)$ for the updated matrix polynomial $Q_\triangle (\lam)$ by  \begin{center}
    $RR_f=\dfrac{\|(M+\triangle M)X_f\Lam_f^2+(D+\triangle D)X_f \Lam_f+(K+\triangle K)X_f\|_F}{\|(M+\triangle M)X_f\Lam_f^2 \|_F+\|(D+\triangle D)X_f \Lam_f \|_F+\|(K+\triangle K)X_f\|_F}$ \\
    
  and  $RR_a=\dfrac{\|(M+\triangle M)X_cP\Lam_a^2+(D+\triangle D)X_cP \Lam_a+(K+\triangle K)X_cP\|_F}{\|(M+\triangle M)X_cP\Lam_a^2 \|_F+\|(D+\triangle D)X_cP \Lam_a \|_F+\|(K+\triangle K)X_cP\|_F}$
\end{center}
respectively \cite{qian2017eigenvalue}. Given $(X_c,\Lam_c)$, we determine the perturbation matrices $\triangle M, \triangle D, \triangle K$ and then calculate the relative residuals to verify the efficiency of reproducing the invariant pairs $(X_cP,\Lam_a)$ and $(X_f,\Lam_f)$ for a structure-preserving perturbed polynomial $Q_\triangle(\lam)$, when the later pair is not known, $\Lam_a$ is given, and $P$ is constructed by the procedure as described in Corollary \ref{Cor:main}.

\begin{example}
Consider the example of a mass-spring system \cite{LiHu} of $10$ degrees of freedom where all the rigid bodies have mass of $1\, Kg$ and all springs have stiffness $1\, kN/m$. Then the quadratic matrix polynomial associated to the model is given by $Q(\lam)=\lam^2 M+\lam D+K$ with $M=I_{10}>0,$
\beano
\small D=\bmatrix{0.4810& -8.3809& 0& 0& 0& 0& 0& 0& 0& 0 \\
-8.3809& 8.3809& -1.0254& 0& 0& 0& 0& 0& 0& 0 \\
0& -1.0254& 1.0254& -7.2827& 0& 0& 0& 0& 0& 0 \\
0& 0& -7.2827& 7.2827& -4.4050& 0& 0& 0& 0& 0 \\
0& 0& 0& -4.4050& 4.4050& -9.9719& 0& 0& 0& 0 \\
0& 0& 0& 0& -9.9719& 9.9719& -5.6247& 0& 0& 0 \\
0& 0& 0& 0& 0& -5.6247& 5.6247& -4.6585& 0& 0 \\
0& 0& 0& 0& 0& 0& -4.6585& 4.6585& -4.1901& 0 \\
0& 0& 0& 0& 0& 0& 0& -4.1901& 4.1901& -2.1160 \\
0& 0& 0& 0& 0& 0& 0& 0& -2.1160& 2.1160}, \eeano 
$$\small K=\bmatrix{2000& -1000& 0& 0& 0& 0& 0& 0& 0& 0 \\
-1000& 3000& -1000& 0& -1000& 0& 0& 0& 0& 0 \\
0& -1000& 2000& -1000& 0& 0& 0& 0& 0& 0 \\
0& 0& -1000& 3000& -1000& 0& 0& -1000& 0& 0 \\
0& -1000& 0& -1000& 3000& -1000& 0& 0& 0& 0 \\
0& 0& 0& 0& -1000& 2000& -1000& 0& 0& 0 \\
0& 0& 0& 0& 0& -1000& 2000& -1000& 0& 0 \\
0& 0& 0& -1000& 0& 0& -1000& 3000& -1000& 0 \\
0& 0& 0& 0& 0& 0& 0& -1000& 2000& -1000 \\
0& 0& 0& 0& 0& 0& 0& 0& -1000& 2000} >0.$$

Computing the eigenvalues of $Q(\lam),$ set $\lam^c_1=-6.7757+71.1468 \mathrm{i},\,\lam^c_2=-6.2938+65.6677 \mathrm{i}.$ Assume $\lam^a_1=-6.16+69.8 \mathrm{i},\,\lam^a_2=-4.7+64.9 \mathrm{i}$. Clearly $\lam^c_1,\,\overline{\lam^c_1},\,\lam^c_2,\,\overline{\lam^c_2}$ are eigenvalues of $Q(\lam)$. Then we want to find the perturbation matrices $\triangle M,\,\triangle D,\,\triangle K$ such that $\lam^a_1,\,\overline{\lam^a_1},\,\lam^a_2,\,\overline{\lam^a_2}$ become eigenvalues of the real symmetric matrix polynomial $Q_{\triangle}(\lam):=\lam^2(M+\triangle M)+\lam(D+\triangle D)+(K+\triangle K),$ whereas the rest of the eigenpairs of $Q_{\triangle}(\lam)$ are same as those unknown (remaining) eigenpairs of $Q(\lam)$ and $\triangle M,\,\triangle K$ must be positive semi-definite. 

Then consider the matrices 
\beano
\Lam_c=\bmatrix{-6.7757&   71.1468&    0&    0\\
  -71.1468&   -6.7757&    0&    0 \\
    0&    0&  -6.2938&   65.6677 \\
    0&    0&  -65.6677&   -6.2938},\,\,\,\,
    \Lam_a=\bmatrix{-6.16& 69.8& 0& 0 \\
-69.8& -6.16& 0& 0 \\
0& 0& -4.7& 64.9 \\
0& 0& -64.9& -4.7} \eeano
\begin{center}
    and \,\,\,\, $X_c=\bmatrix{-0.28211&  -0.08966&  -0.29723&   0.07950 \\
   0.83728&   0.03582&   0.60914&  -0.29543 \\
  -0.59676&  -0.02745&  -0.10053&   0.21115 \\
   1.00000&  0&  -0.26137&  -0.02870 \\
  -0.91905&  -0.05869&  -0.44955&   0.10855 \\
   0.21414&   0.23271&   0.45059&   0.21796 \\
   0.16581&  -0.13691&  -0.60870&  -0.17516 \\
  -0.63899&   0.23427&   1.00000&   0 \\
   0.23640&  -0.07021&  -0.51521&   0.00617 \\
  -0.07220&   0.03048&   0.21217&  -0.03517},$ \end{center}
  such that $(X_c, \Lam_c)$ is an invariant pair of $Q(\lam).$
  
By Remark \ref{Reco_Chu_postdefn} we have $\gamma_1=4.2986 - 485.9606 \mathrm{i},\,\gamma_2=20.523 - 319.028 \mathrm{i},\,\alpha_1=2.1493,\,\alpha_2=10.2615,\,\beta_1=242.9803,\,\beta_2=159.514$ and taking $a_1=0.00098,\,a_2=0.00831$ we have  $$P=\bmatrix{0.00098&  -0.96441&   0& 0\\
   1.18488&  -0.00098&   0& 0\\
   0& 0&   0.00831&  -5.00141 \\
  0& 0&   5.25988&  -0.00831}.$$ Consequently, we obtain the real symmetric positive definite matrix $$Z=\bmatrix{
  0.07448&  -0.00328&  0.00042&  -0.00065\\
  -0.00328&  0.06075&  0.00080&  -0.00029\\
  0.00042&  0.00080&  0.02586&  -0.01101\\
  -0.00065&  -0.00029&  -0.01101&  0.01727}.$$ Hence on applying Remark \ref{Reco_Chu_postdefn} we obtain the real symmetric matrices as
\beano
\scriptsize \triangle M=10^{-3} \bmatrix{
9.3095&  -24.4279&   14.5031&  -18.7431&   23.9416&   -8.4718&    1.6675&    3.4774&   -0.0641&   -0.6650 \\
  -24.4279&   67.9437&  -41.9671&   57.9737&  -67.6659&   19.4633&    0.7605&  -20.6713&    4.7560&    0.0216 \\
   14.5031&  -41.9671&   28.1813&  -43.3236&   43.8138&  -10.7203&   -5.1443&   23.2965&   -7.8302&    1.9365 \\
  -18.7431&   57.9737&  -43.3236&   75.9119&  -65.2624&   12.7806&   16.1761&  -54.3748&   20.9451&   -6.8164 \\
   23.9416&  -67.6659&   43.8138&  -65.2624&   69.7849&  -19.2399&   -4.4682&   30.3068&   -9.2572&    1.8447 \\
   -8.4718&   19.4633&  -10.7203&   12.7806&  -19.2399&   10.4451&   -4.9098&    3.0152&   -2.2823&    1.3902 \\
    1.6675&    0.7605&  -5.1443&   16.1761&   -4.4682&   -4.9098&   11.1734&  -24.1222&   10.8057&   -4.2638 \\
    3.4774&  -20.6713&   23.2965&  -54.3748&   30.3068&    3.0152&  -24.1222&   60.4317&  -25.8859&    9.8290\\
   -0.0641&    4.7560&   -7.8302&   20.9451&   -9.2572&   -2.2823&   10.8057&  -25.8859&   11.4607&   -4.4693\\
   -0.6650&    0.0216&    1.9365&   -6.8164&    1.8447&    1.3902&   -4.2638&    9.8290&   -4.4693&    1.8039} \geq 0, \eeano
$$\scriptsize \triangle D=\bmatrix{
0.1225&  -0.2216&   0.2315&  -0.4845&   0.1814&  -0.0467&  -0.2684&   0.4364&  -0.2122&   0.0689\\
  -0.2216&   0.1724&  -0.3101&   0.7276&  -0.0057&   0.0788&   0.5863&  -0.6878&   0.4014& -0.1192\\
   0.2315&  -0.3101&   0.2840&  -0.5046&   0.2195&  -0.2047&  -0.2389&   0.2513&  -0.1695&   0.0314\\
  -0.4845&   0.7276&  -0.5046&   0.7066&  -0.6532&   0.6396&   0.0267&   0.0979&   0.0462&   0.0496\\
   0.1814&  -0.0057&   0.2195&  -0.6532&  -0.0649&  -0.1825&  -0.5081&   0.6043&  -0.3935&   0.1195\\
  -0.0467&   0.0788&  -0.2047&   0.6396&  -0.1825&   0.0679&   0.3234&  -0.7218&   0.3438&  -0.1373\\
  -0.2684&   0.5863&  -0.2389&   0.0267&  -0.5081&   0.3234&  -0.3182&   0.5763&  -0.2716&   0.1421\\
   0.4364&  -0.6878&   0.2513&   0.0979&   0.6043&  -0.7218&   0.5763&  -1.1491&   0.4605&  -0.2451\\
  -0.2122&   0.4014&  -0.1695&   0.0462&  -0.3935&   0.3438&  -0.2716&   0.4605&  -0.1955&   0.1051\\
   0.0689&  -0.1192&   0.0314&   0.0496&   0.1195&  -0.1373&   0.1421&  -0.2451&   0.1051&  -0.0505},$$ 
$$\scriptsize \triangle K=\bmatrix{
28.0577&   -81.5942&    52.7467&   -81.4589&    88.5228&   -28.0835&    -1.6524&    38.1178&   -11.1266&     2.8698 \\
   -81.5942&   241.9383&  -160.7575&   255.6535&  -262.6384&    75.0079&    19.1931&  -137.2175&    43.9081&   -12.5844 \\
    52.7467&  -160.7575&   112.0746&  -186.2415&   175.5596&   -42.0997&   -27.4937&   117.4612&   -41.1852&    12.9837 \\
   -81.4589&   255.6535& -186.2415&   323.4317&  -283.4432&    57.5522&    64.8610&  -228.1801&    84.4493&   -27.9849 \\
    88.5228&  -262.6384&   175.5596&  -283.4432&   290.1658&   -85.8764&   -19.0487&   152.9206&   -49.7268&    14.6872 \\
   -28.0835&    75.0079&   -42.0997&    57.5522&   -85.8764&    42.9652&   -22.8915&    -0.1737&    -5.3397&     3.0183 \\
    -1.6524&    19.1931&   -27.4937&    64.8610&   -19.0487&   -22.8915&    49.3152&   -88.8186&    38.7107&   -14.2052 \\
    38.1178&  -137.2175&   117.4612&  -228.1801&   152.9206&    -0.1737&   -88.8186&   214.3873&   -87.0151&    30.7811 \\
   -11.1266&    43.9081&   -41.1852&    84.4493&   -49.7268&    -5.3397&    38.7107&   -87.0151&    36.2412&   -13.0667 \\
     2.8698&   -12.5844&    12.9837&   -27.9849&    14.6872&     3.0183&   -14.2052&    30.7811&   -13.0667&     4.7860} \geq 0.$$  
      
Let $(X_f,\Lam_f)$ denote the `fixed' (unknown) invariant pair of $Q(\lam).$ Then the relative residuals of $(X_f,\Lam_f),\,(X_cP,\Lam_a)$ for the updated system $Q_\triangle (\lam)$ are given by 
$$RR_f=7.5511 \times 10^{-14} \approx 0,\,\,\,\,RR_a=4.6172 \times 10^{-14} \approx 0,$$  which ensures that $\lam^a_1,\,\overline{\lam^a_1},\,\lam^a_2,\,\overline{\lam^a_2}$ are the eigenvalues of $Q_{\triangle}(\lam)$ and rest of the eigenpairs of $Q_{\triangle}(\lam)$ are same as those unknown eigenpairs of $Q(\lam)$. 

In the following figures we plot the relative residuals of $(X_f,\Lam_f)$ and $(X_cP,\Lam_a)$ for the updated system $Q_\triangle(\lam)$ for different parametric values of $a_1,\, a_2$ as mentioned in Corollary \ref{Cor:main} for the case when $Q(\lam) \in \Q_n(T,1,1) \subset \R^{n \times n}[\lam]$. We plot $RR_f,\,RR_a$ by choosing $a_1=0.00098,\,a_2=0.01(j-5)$ and $a_1=-0.0024,\,a_2=0.01(j-5),\,j=1,\hdots,10$ in Figure \ref{fig1}.  

\begin{figure}[H] 
\centering
\subfigure[Relative residuals corresponding to the fixed eigenpairs]{\includegraphics[height=6 cm,width=6.5 cm]{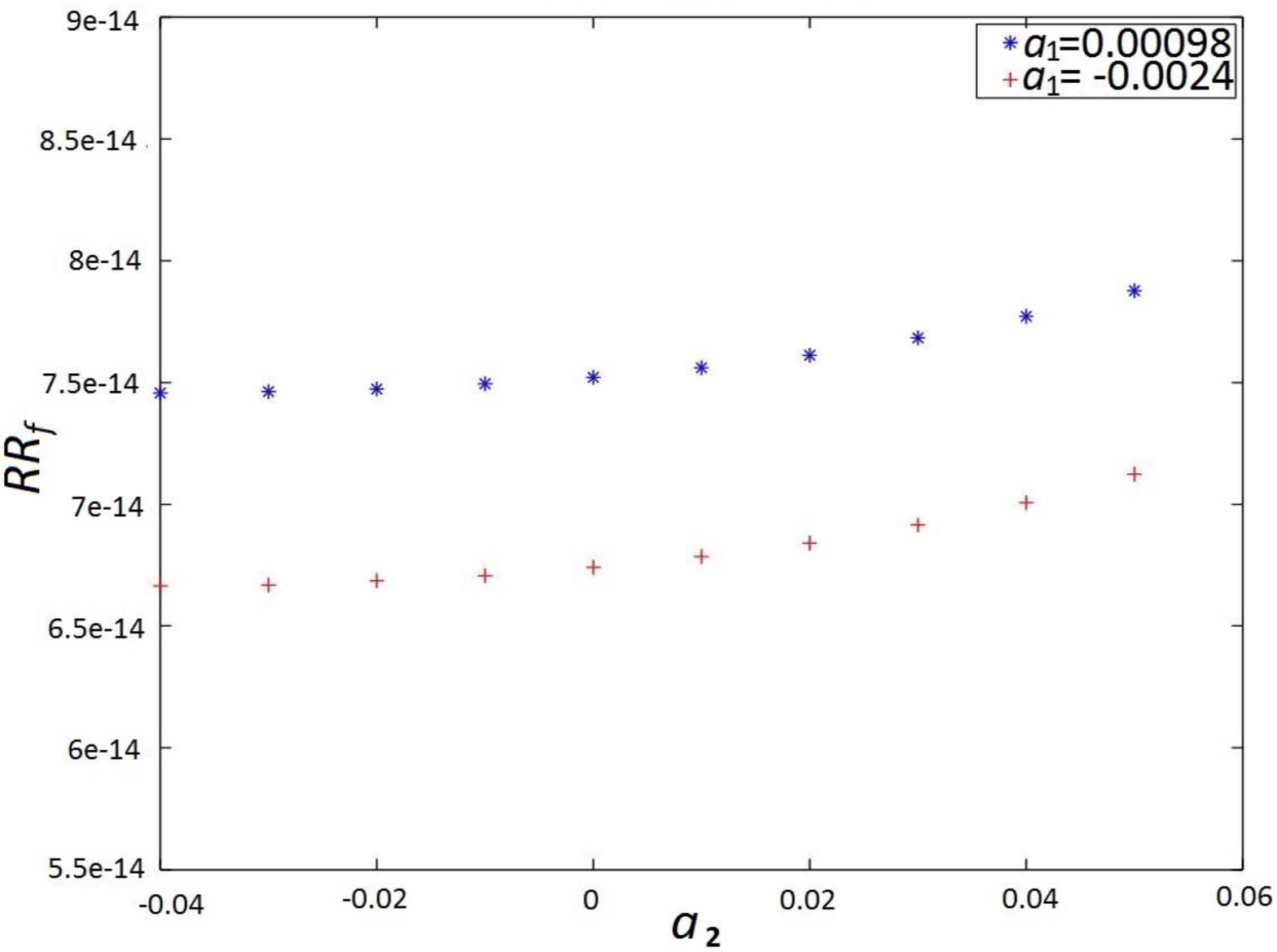}} \hspace{0.3cm}
\subfigure[Relative residuals corresponding to the aimed eigenpairs]{\includegraphics[height=6 cm,width=6.5 cm]{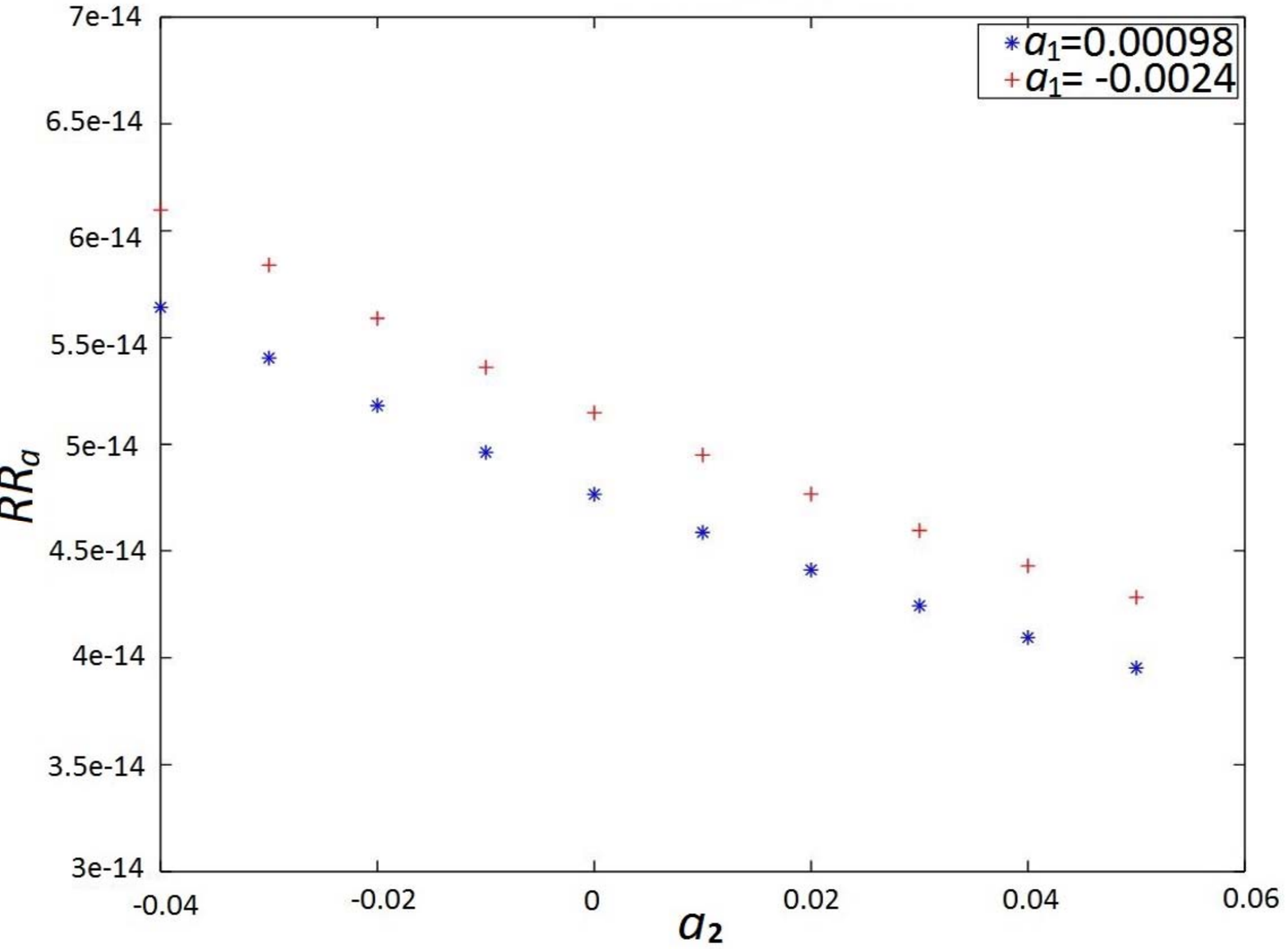}}
\caption{Relative residuals for different perturbations corresponding to the parameters $a_1, a_2$ which define the parameter matrix $P$ for the perturbations $\triangle M, \triangle D, \triangle K$.} \label{fig1}
\end{figure} 

  
\end{example}
\begin{example}
We consider this example from \cite{Mao_Dai}. Consider the $T-$even quadratic matrix polynomial $Q(\lam)=\lam^2 M+\lam D+K$ with 
\beano
M=I_{3}>0,\,\,
D=\bmatrix{0& -2& 4 \\ 2& 0& -2 \\ -4& 2& 0}, \,\,
K=\bmatrix{13& 2& 1 \\ 2& 7& 2 \\ 1& 2& 4} >0.
\eeano 
Let $\lam^c_1=0.8878 \mathrm{i},\,\lam^c_2=3.1895 \mathrm{i}$ and $\lam^a_1=2 \mathrm{i},\,\lam^a_2=3.5 \mathrm{i}$. Then it is easy to verify that $\lam^c_1,\,\overline{\lam^c_1},\,\lam^c_2,\,\overline{\lam^c_2}$ are eigenvalues of $Q(\lam)$. Now it is necessary to determine the real perturbation matrices $\triangle M,\,\triangle D,\,\triangle K$ in such a manner that $\lam^a_1,\,\overline{\lam^a_1},\,\lam^a_2,\,\overline{\lam^a_2}$ become the eigenvalues of the $T-even$ quadratic matrix polynomial $Q_{\triangle}(\lam):=\lam^2(M+\triangle M)+\lam(D+\triangle D)+(K+\triangle K),$ whereas rest of the eigenpairs of $Q_{\triangle}(\lam)$ are same as those unknown eigenpairs of $Q(\lam)$. Then applying Remark \ref{Reco_Mao_Dai} we have
\beano
\Lam_c=\bmatrix{0&   0.8878&   0& 0 \\
  -0.8878&   0&  0& 0 \\
   0& 0&   0&  3.1895 \\
   0& 0&  -3.1895&   0},\,\,\,\,
    \Lam_a=\bmatrix{0&   2&   0& 0\\
  -2&   0&   0&   0 \\
   0&   0& 0&   3.5 \\
   0&  0&  -3.5&   0} \eeano
\begin{center}
    and \,\,\,\, $X_c=\bmatrix{-0.07809&  -0.42447&   0.36993&   0.04128 \\
  -0.41817&   0.44484&   1&   0 \\
   1&  0&   0.46926&   0.27553}$ \end{center}
 such that $(X_c,\Lam_c)$ is an invariant pair of $Q(\lam)$.
Using Remark \ref{Reco_Mao_Dai} we choose the matrix $$P=\bmatrix{
-0.64146&  -0.87909&   0& 0\\
   1.27500&  -1.27500&   0& 0\\
   0& 0&   0.55689&   0.99159 \\
   0& 0&  -0.90016&   0.90016}$$ and obtain the real symmetric matrix $$Z=\bmatrix{
  -0.06974&  -0.05987&   0.00391&   0.00614 \\
  -0.05987&  -0.46773&   0.00287&   0.02937 \\
   0.00391&   0.00287&  -0.00892&  -0.20898 \\
   0.00614&   0.02937&  -0.20898&  -0.12565}.$$ 
   
   Therefore, we obtain the real matrices $\triangle M,\,\triangle D,\,\triangle K$ so that 
\beano
M+\triangle M=0.1 \bmatrix{
9.0132&  0.6435&  -0.0005 \\
  0.6435&  9.0789&  -0.5260 \\
  -0.0005&  -0.5260&  8.7177},\, 
D+\triangle D=\bmatrix{
0&  -1.6534&  2.8543\\
  1.6534&  0&  -1.2119 \\
  -2.8543&  1.2119&  0},\eeano 
$$K+\triangle K=\bmatrix{
15.1304&    2.0337&    0.2838 \\
    2.0337&    9.5026&   -0.3897 \\
    0.2838&   -0.3897&    9.2174}$$  
and $\|\triangle M\|_F=0.2202,\,\|\triangle D\|_F=2.0268,\,\|\triangle K\|_F=7.1044$.       
Let $(X_f,\Lam_f)$ denote the invariant pair corresponding to the fixed eigenpairs of $Q(\lam).$ Then the relative residuals of $(X_f,\Lam_f),\,(X_cP,\Lam_a)$ for the updated polynomial $Q_\triangle (\lam)$ are given by 
$$RR_f=7.0842 \times 10^{-16} \approx 0,\,\,\,\,RR_a=3.2018 \times 10^{-16} \approx 0,$$ which ensures that the unknown eigenpairs of $Q(\lam)$ remains to be the eigenpairs of $Q_{\triangle}(\lam)$ and $\lam^a_1,\,\overline{\lam^a_1},\,\lam^a_2,\,\overline{\lam^a_2}$ are eigenvalues of the $T-$even matrix polynomial $Q_{\triangle}(\lam)$. Hence, eigenvalues are replaced successfully with maintaining no spillover on unmeasured eigenpairs of $Q(\lam)$.  

Moreover, choosing $$P=\bmatrix{
-0.64146&  -0.87909&   0& 0\\
   -c_1&  c_1&   0& 0\\
   0& 0&   0.55689&   0.99159 \\
   0& 0&  -c_2&   c_2}$$
   we plot the relative residuals of $(X_f,\Lam_f)$ and $(X_cP,\Lam_a)$ for the updated matrix polynomial $Q_{\triangle}(\lam)$ for various parametric values of $c_1,\,c_2$ as mentioned in Remark \ref{Reco_Mao_Dai}. In Figure \ref{fig2}, we plot $RR_f,\,RR_a$ choosing $c_1=-1.2750,\,c_2=0.2j$ and $c_1=11.0214,\,c_2=0.2j,\,j=1,\hdots,10$.

\begin{figure}[H] 
\centering
\subfigure[Relative residuals corresponding to the fixed eigenpairs]{\includegraphics[height=6 cm,width=6.5 cm]{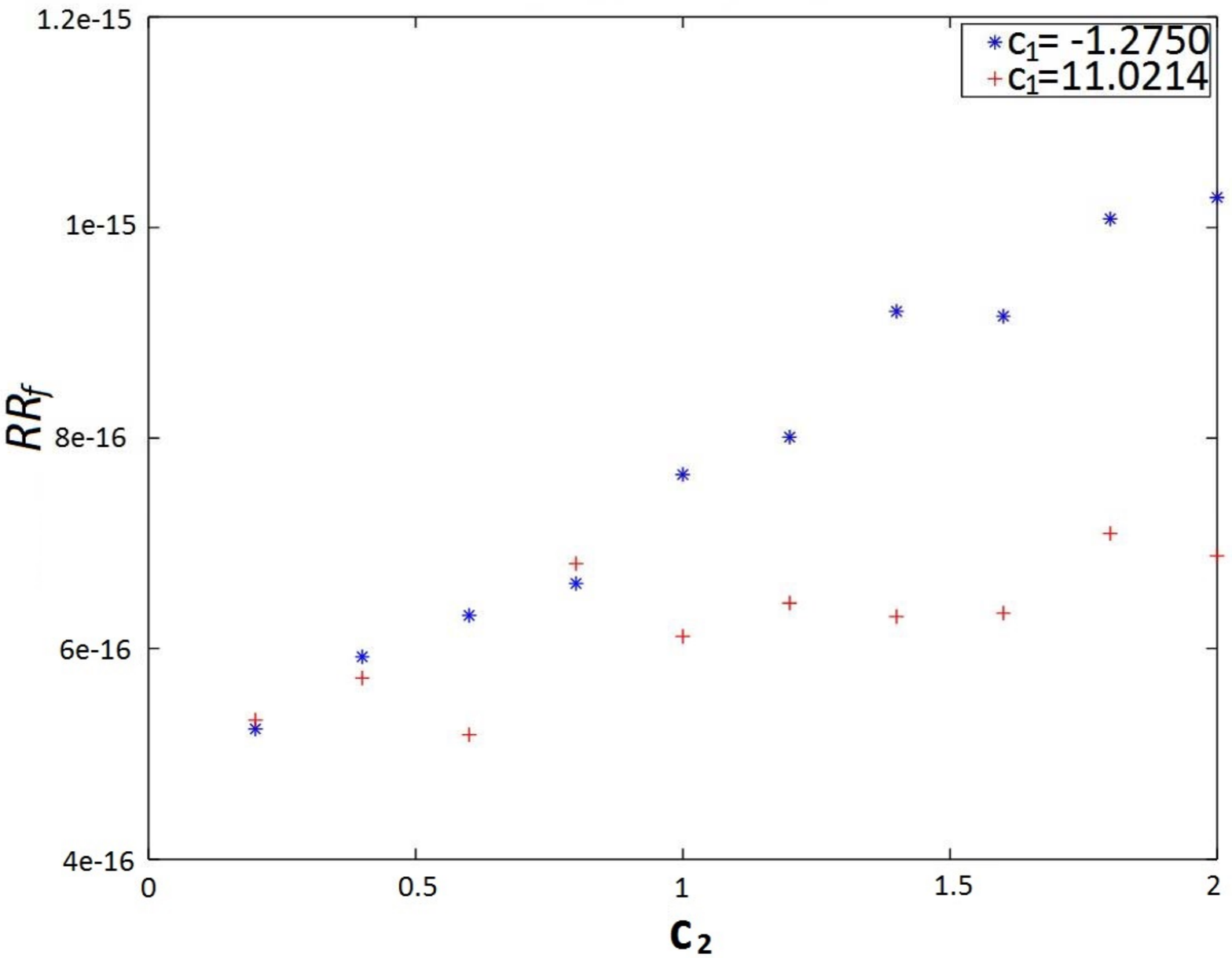}} \hspace{0.3cm}
\subfigure[Relative residuals corresponding to the aimed eigenpairs]{\includegraphics[height=6 cm,width=6.5 cm]{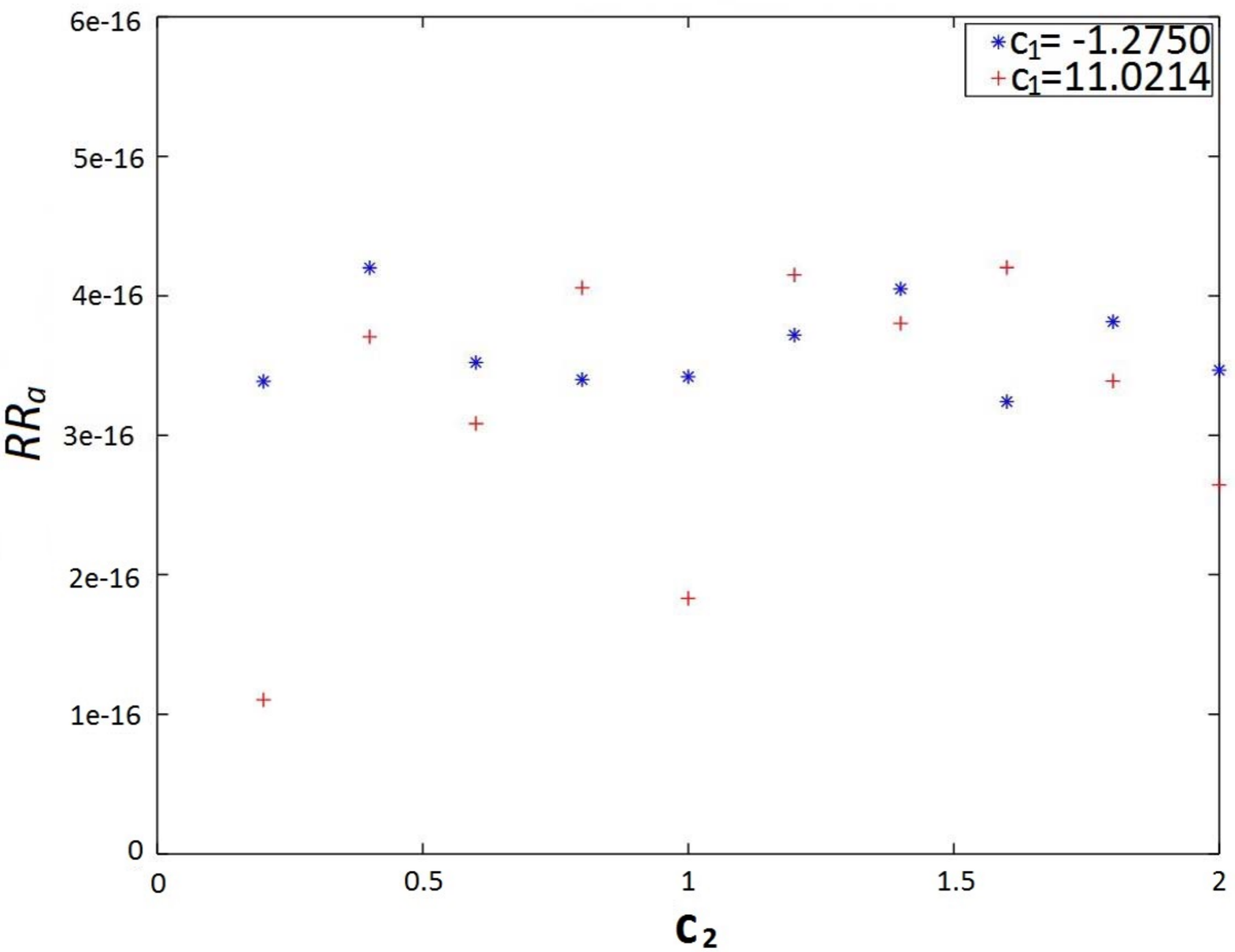}}
\caption{Relative residuals for different perturbations corresponding to the parameters $c_1, c_2$ which define the parameter matrix $P$ for the perturbations $\triangle M, \triangle D, \triangle K$.} \label{fig2}
\end{figure}

  
\end{example}

\begin{example}
In this example, we consider the matrices $M,\,K$ as randomly generated real symmetric positive definite matrices of order $52 \times 52$ and $D$ is a randomly generated real skew-symmetric matrix of order $52 \times 52$. Let $\lam^c_1=46.76551\mathrm{i},\,\lam^c_2=16.58514 \mathrm{i},\,\lam^c_3=14.33130 \mathrm{i},\,\lam^c_4=8.44632 \mathrm{i}$ and $\lam^a_1=2.1145 \mathrm{i},\,\lam^a_2=0.2374 \mathrm{i},\,\lam^a_3=11.3266 \mathrm{i},\,\lam^a_4=0.0752 \mathrm{i}$. Suppose we want to replace the $8$ eigenvalues $\lam^c_1,\,-\lam^c_1,\,\lam^c_2,\,-\lam^c_2,\,\lam^c_3,\,-\lam^c_3,\,\lam^c_4,\,-\lam^c_4$ of $Q(\lam)=\lam^2 M+\lam D+K$ by the desired scalars $\lam^a_1,\,-\lam^a_1,\,\lam^a_2,\,-\lam^a_2,\,\lam^a_3,\,-\lam^a_3,\,\lam^a_4,\,-\lam^a_4$ with maintaining no spillover effect on unmeasured eigenpairs of $Q(\lam)$. Then applying Remark \ref{Reco_Mao_Dai} we define the matrices $\Lam_c,\,\Lam_a,\,X_c$ accordingly and choosing the matrix $$P=\diag \left(\bmatrix{0.24821&   0.10510\\ -0.00973&   0.00973},\,\bmatrix{-0.13107&   0.15040\\ -0.04713&   0.04713},\,\bmatrix{0.11182&   0.05132\\-0.10860&   0.10860},\,\bmatrix{0.00660&  -0.05718\\ 0.16688&  -0.16688} \right)$$ we obtain the real perturbation matrices $\triangle M,\,\triangle D,\,\triangle K$ with $\|\triangle M\|_F=6.1883,\,\|\triangle D\|_F=52.153,\linebreak \|\triangle K\|_F= 189.74$ and $\|\triangle M-\triangle M^T\|_F=7.4674\times 10^{-13},\,\|\triangle D+\triangle D^T\|_F=1.4374 \times 10^{-11},\,\|\triangle K-\triangle K^T\|_F=2.2838 \times 10^{-11}$, which ensures that $\triangle M,\,\triangle K$ are symmetric matrices while $\triangle D$ is a skew-symmetric matrix. It should be noted that $RR_a=1.1668 \times 10^{-12}$ is nearly zero, which implies that $\lam^a_1,\,-\lam^a_1,\,\lam^a_2,\,-\lam^a_2,\,\lam^a_3,\,-\lam^a_3,\,\lam^a_4,\,-\lam^a_4$ are eigenvalues of $Q_\triangle (\lam)=\lam^2 (M+\triangle M)+\lam (D+\triangle D)+(K+\triangle K)$. Let $(X_f,\Lam_f)$ denote the fixed invariant pair corresponding to the fixed eigenpairs of $Q(\lam).$ Then we have  
$$\dfrac{\|MX_f\Lam_f^2+D X_f \Lam_f+K X_f\|_F}{\|MX_f\Lam_f^2 \|_F+\|D X_f \Lam_f \|_F+\|KX_f\|_F}= 5.0119 \times 10^{-14}.$$ Besides, 
$RR_f=3.1912 \times 10^{-13}$ is nearly zero, which guarantees that fixed eigenpairs of $Q(\lam)$ remain the eigenpairs of $Q_\triangle (\lam)$. Hence, it shows that a few eigenvalues of $Q(\lam)$ are replaced by the desired scalars with maintaining no spillover. 

However, choosing $$P=\diag \left(\bmatrix{0.24821&   0.10510\\ -0.00973&   0.00973},\,\bmatrix{-0.13107&   0.15040\\ -0.04713&   0.04713},\,\bmatrix{0.11182&   0.05132\\-0.10860&   0.10860},\,\bmatrix{a_4&  -0.05718\\ -c_4&   c_4} \right)$$ the relative residuals of $(X_f,\Lam_f)$ and $(X_cP,\Lam_a)$ for the updated system $Q_\triangle (\lam)$ has been plotted in the following figures for several parametric values of $a_4,\,c_4$ as given in Remark \ref{Reco_Mao_Dai}. In particular, we plot  $RR_f,\,RR_a$ choosing $c_4=6.52,\,a_4=0.02j$ and $c_4=0.2654,\,a_4=0.02j,\,j=1,\hdots,10$ in Figure \ref{fig3}.  

\begin{figure}[H] 
\centering
\subfigure[Relative residuals corresponding to the fixed eigenpairs]{\includegraphics[height=6 cm,width=6.5 cm]{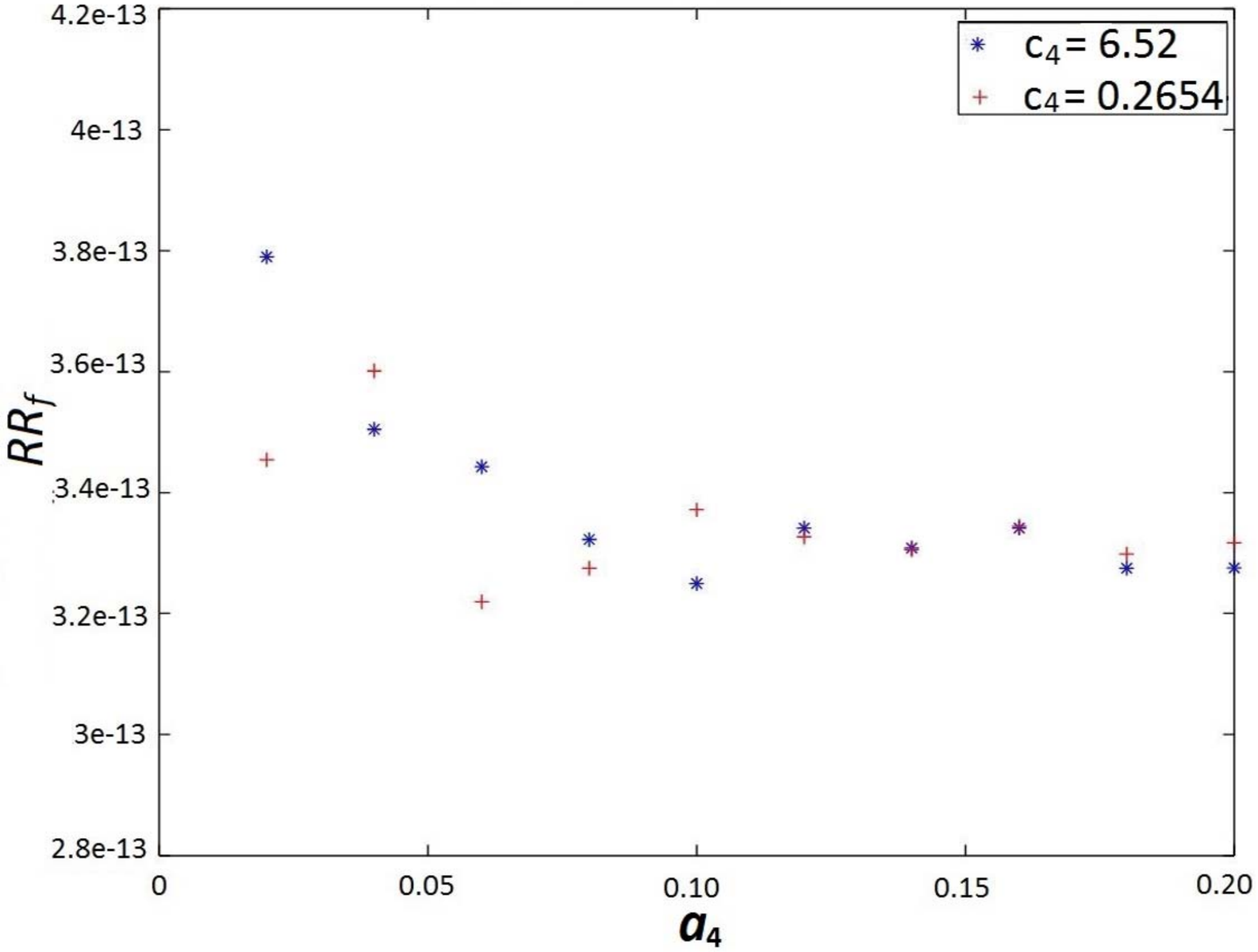}} \hspace{0.3cm}
\subfigure[Relative residuals corresponding to the aimed eigenpairs]{\includegraphics[height=6 cm,width=6.5 cm]{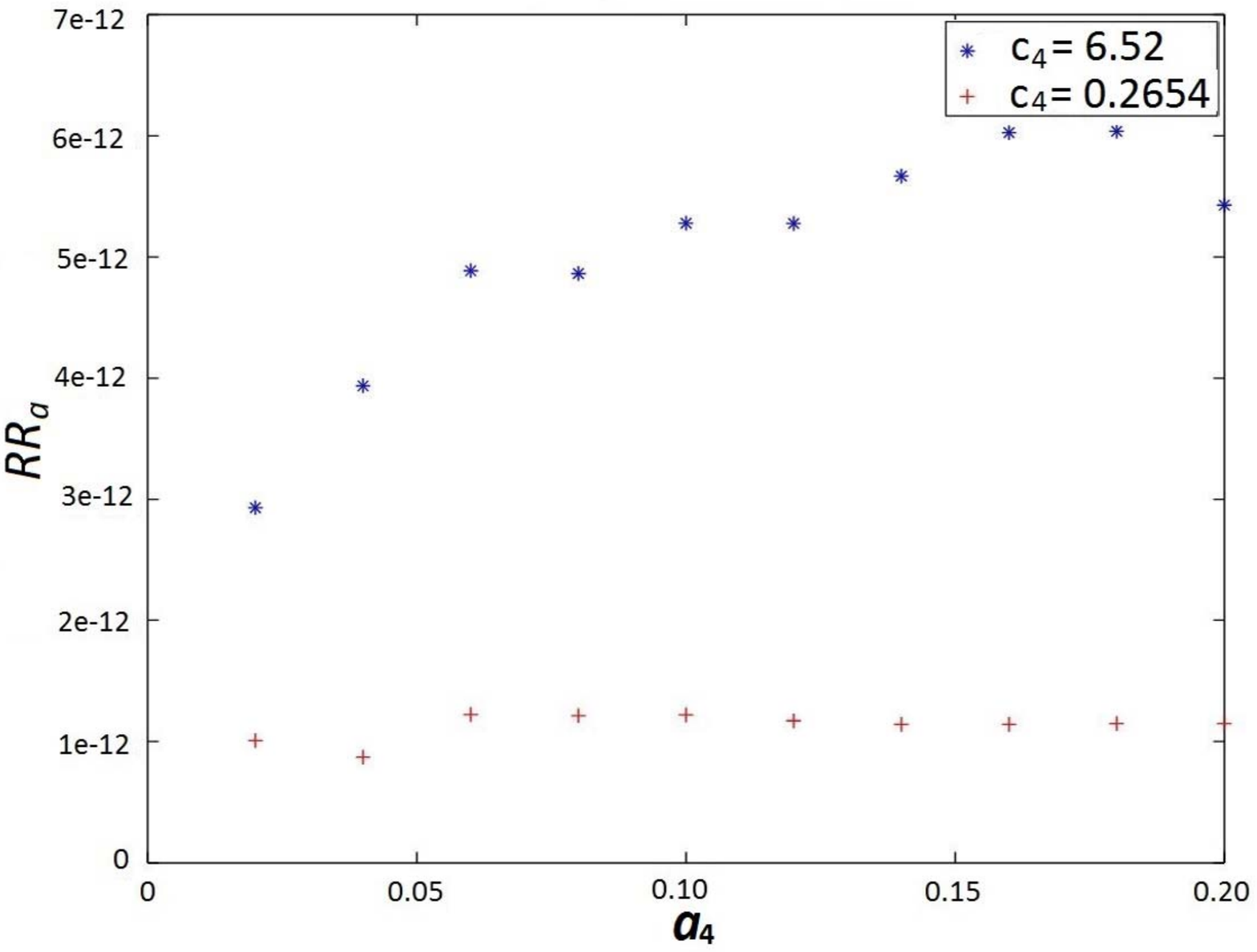}}
\caption{Relative residuals for different perturbations corresponding to the parameters $a_4, c_4$ which define the parameter matrix $P$ for the perturbations $\triangle M, \triangle D, \triangle K$.} \label{fig3}
\end{figure} 

\end{example}

{\bf Conclusion.} We consider the \textit{structure-preserving eigenvalue embedding problem} (SEEP) for regular quadratic polynomials with symmetry structures. First we derive structure-preserving perturbations of a structured qudratic polynomial that reproduce a desired invariant pair and preserve an invariant pair (need not be known) of the unperturbed polynomial. Then we utilize these results for solving the SEEP for quadratic structured matrix polynomials which include symmetric, Hermitian, $\star$-odd and $\star$-even matrix polynomials. We show that the obtained solutions for SEEP correspond to existing results in the literature for certain structured matrix polynomials that arise in real-world applications. Finally, we illustrate the applicability of the obtained results through numerical examples.  


\end{document}